\documentclass[a4paper,12pt,leqno, twoside]{article}
\usepackage{amssymb}
\usepackage{amsthm}
\usepackage{amstext}
\usepackage{amsbsy}
\usepackage{amscd}
\usepackage{amsopn}
\usepackage{amsmath}
\usepackage{mathtools}
\usepackage{mathdots}
\usepackage[mathscr]{eucal}
\usepackage{graphicx}
\usepackage{epic,eepic}
\usepackage{epsfig}
\usepackage[all]{xy}
\usepackage[french]{babel}
\usepackage[latin1]{inputenc}

    \def\AM{{\mathbb{A}}}

\def\DG{{\mathfrak D}}    
    
    \def\FM{{\mathbb{F}}}
    \def\GM{{\mathbb{G}}}

    \def\NM{{\mathbb{N}}}
    
    \def\PM{{\mathbb{P}}}
    \def\QM{{\mathbb{Q}}}
    
\def\SG{{\mathfrak S}}

    \def\ZM{{\mathbb{Z}}}

\def\Bb{{\mathbf B}}

    \def\FC{{\mathcal{F}}}
\def\Gb{{\mathbf G}}  \def\gb{{\mathbf g}}

\def\Lb{{\mathbf L}}    
  \def\mb{{\mathbf m}}  \def\MC{{\mathcal{M}}}

\def\Pb{{\mathbf P}}    
\def\Qb{{\mathbf Q}}

\def\Tb{{\mathbf T}}    
\def\Ub{{\mathbf U}}    \def\UC{{\mathcal{U}}}
\def\Vb{{\mathbf V}}

          \def\cdo{{\dot{c}}}

          \def\wdo{{\dot{w}}}

          \def\fba{{\bar{f}}}

          \def\pba{{\bar{p}}}

\def\a{\alpha}
\def\b{\beta}

\def\G{\Gamma}
\def\d{\delta}
\def\D{\Delta}

\def\l{\lambda}
\def\L{\Lambda}

\def\o{\omega}
\def\O{\Omega}

\def\x{\xi}

\def\xib{{\boldsymbol{\xi}}}

\DeclareMathOperator{\Hom}{{\mathrm{Hom}}}

\DeclareMathOperator{\Ker}{{\mathrm{Ker}}}

\DeclareMathOperator{\res}{{\mathrm{res}}}

\DeclareMathOperator{\rg}{{\mathrm{rg}}}

\DeclareMathOperator{\Spec}{{\mathrm{Spec}}}

\def\oFq{\o{\FM}_q}

\def\ssi{si et seulement si }
\def\para{sous-groupe parabolique }
\def\paras{sous-groupes paraboliques }

\def\borel{sous-groupe de Borel }
\def\borels{sous-groupes de Borel }

\def\Omeq{\Omega_{\FM_q}}

\def\opp{\rm opp}
\def\norm{\mathrm {Norm}}

\def\ord{{\rm ord}}

\def\re{{\rm red}}

\def\simto{\buildrel\hbox{$\sim$}\over\longrightarrow}
\def\vide{\varnothing}
\def\leq{\leqslant}
\def\geq{\geqslant}
\def\into{\hookrightarrow}
\def\onto{\twoheadrightarrow}
\def\id{\mathop{\mathrm{Id}}\nolimits}

\def\ba{\backslash}
\def\wt{\widetilde}

\def\o#1{\overline{#1}}

\def\To#1{\buildrel\hbox{\tiny{$#1$}}\over\longrightarrow}

\def\to{\rightarrow}

\def\Hom{\mathop{\hbox{\rm Hom}}\nolimits}

\def\dim{\mathop{\mbox{\rm dim}}\nolimits}

\def\ini{\setcounter{equation}{\value{subsubsection}}\addtocounter{subsubsection}{1}}

\makeatletter
\renewcommand{\subsubsection}{\@startsection{subsubsection}{3}{\parindent}{-\baselineskip}{-0.01\baselineskip}{\bf}}
\renewcommand*{\@seccntformat}[1]{%
  \csname the#1\endcsname\
} \makeatother

\def\ali{\subsubsection{}\setcounter{equation}{0}}

\newtheoremstyle{th}
  {\baselineskip}{.5\baselineskip}{\itshape}
  {\parindent}{\bf}
  {--}{.5em}{}

\newtheoremstyle{def}
  {\baselineskip}{\baselineskip}{}
  {\parindent}{\bf}
  {--}{.5em}{}

\newtheoremstyle{th*}
  {.5\baselineskip}{.5\baselineskip}{\itshape}
  {\parindent}{\bf}
  {--}{.5em}{}

\newtheoremstyle{remark*}
  {.5\baselineskip}{.5\baselineskip}{}
  {\parindent}{\bf}
  {--}{.5em}{}

\newtheoremstyle{remark}
  {.5\baselineskip}{.5\baselineskip}{}
  {\parindent}{\bf}
  {--}{.5em}{}

\swapnumbers \theoremstyle{th}
\newtheorem{theo}[subsubsection]{\sc Th{\'e}or{\`e}me.\bf}
\newtheorem{lemme}[subsubsection]{\sc Lemme.\bf}
\newtheorem{prop}[subsubsection]{\sc Proposition.\bf}
\newtheorem{coro}[subsubsection]{\sc Corollaire.\bf}

\theoremstyle{def}

\theoremstyle{remark}

\theoremstyle{th*}
\newtheorem*{thm}{\sc Th{\'e}or{\`e}me.}
\newtheorem*{lem}{\sc Lemme.}

\theoremstyle{remark*}
\newtheorem*{rem}{\sc Remarque.}

\newcommand{\findem}{\hfill$\Box$\par\medskip}
\newcommand{\dem}{\indent {\it Preuve :} \rm }

\newenvironment{preuve}{\dem}{\findem}

\title{Sur la cohomologie de la compactification des vari\'et\'es de Deligne-Lusztig}
\author{Haoran Wang}
\date{}

\theoremstyle{plain}

\begin{document}
\maketitle

\renewcommand{\proofname}{\indent Preuve}

\def\dd{D_d^\times}
\def\mdro{\MC_{Dr,0}}
\def\mdrn{\MC_{Dr,n}}
\def\mdr{\MC_{Dr}}
\def\mlto{\MC_{LT,0}}
\def\mltn{\MC_{LT,n}}
\def\mlt{\MC_{LT}}
\def\mltK{\MC_{LT,K}}
\def\LJ{{\rm LJ}}
\def\JL{{\rm JL}}
\def\SL{{\rm SL}}
\def\GL{{\rm GL}}
\def\Spf{{\rm Spf}}
\def\DL{{\rm DL}}
\def\PGL{{\rm PGL}}
\def\Coef{{\rm Coef}}

\def\Ql{\QM_{\ell}}
\def\Zl{\ZM_{\ell}}
\def\Fl{\FM_{\ell}}
\abstract{
Nous \'etudions la cohomologie de la compactification des vari\'et\'es de Deligne-Lusztig associ\'ees aux \'el\'ements de Coxeter. Nous pr\'esentons une conjecture des relations entre la cohomologie de la vari\'et\'e et la cohomologie de ses compactifications partielles. Nous prouvons la conjecture dans le cas du groupe lin\'eaire g\'en\'eral.}

\tableofcontents
\section{Introduction}
Soient $q$ une puissance d'un nombre premier $p,$ $\Gb$ un groupe r\'eductif connexe sur $\o{\FM}_p$ muni d'une structure $\FM_q$-rationnelle et $F$ l'isog\'enie de Frobenius correspondante. Si $w$ est un \'el\'ement du groupe de Weyl, dans leur article fondateur \cite{deligne-lusztig} Deligne et Lusztig ont introduit deux vari\'et\'es $X(w)$ et $Y(w)$ munies d'une action de $\Gb^F$ ainsi qu'un morphisme fini \'etale $\Gb^F$-\'equivariant $\pi:Y(w)\to X(w),$ faisant de $X(w)$ un quotient de $Y(w)$ par l'action d'un groupe fini commutatif (dans \cite{deligne-lusztig} la vari\'et\'e $Y(w)$ est not\'ee $\wt{X}(\wdo)$). Dans leur travail, ils ont \'egalement construit une compactification lisse $\o{X}(w)$ de $X(w)$ \`a la Bott-Samelson-Demazure-Hansen.

Lorsque $w$ est un \'el\'ement de Coxeter, l'\'etude de la g\'eom\'etrie des vari\'et\'es de Deligne-Lusztig a commenc\'e par Lusztig \cite{lusztig-coxeter}, et puis Bonnaf\'e et Rouquier \cite{bonnafe-rouquier-coxeter}, et Dudas \cite{dudas}. Dans ce cas, la compactification $\o X(w)$ poss\`ede une stratification dont les strates sont index\'ees par les \paras propres $F$-stables. Plus pr\'ecis\'ement, soient $\Pb$ un \para propre $F$-stable et $\Ub$ son radical unipotent, la strate ferm\'ee associ\'ee \`a $\Pb$ est
$$
\o{X}_{\Pb}(w):=\o X(w)^{\Ub^F}.
$$
Notons $j:X(w)\into \o X(w)$ et $i_\Pb:\o X_\Pb(w)\into \o X(w)$ les immersions naturelles.

Le but principal de cet article est de d\'emontrer le th\'eor\`eme suivant (voir \ref{main theorem 2}):
\begin{thm}
Soient $\Gb=\GL_d,$ $F$ l'endomorphisme de Frobenius standard $(a_{ij})_{1\leq i,j\leq d}\mapsto (a_{ij}^q)_{1\leq i,j\leq d}$ et $w=(1,\ldots,d)\in\SG_d.$ Alors le morphisme de restriction:
$$
R\G(X(w),\pi_*\L)=R\G(\o X(w) ,Rj_*(\pi_*\L))\To{\res.}R\G(\o{X}_{\Pb}(w),i^*_\Pb Rj_*(\pi_*\L))
$$
induit un isomorphisme
\begin{equation}\label{hh}
R\G(X(w),\pi_*\L)^{\Ub^F}\simto R\G(\o{X}_{\Pb}(w),i^*_\Pb Rj_*(\pi_*\L)),
\end{equation}
o\`u $\L=\ZM/\ell^m$ pour un nombre premier $\ell\neq p.$
\end{thm}
Notre motivation pour ce th\'eor\`eme vient du lien avec les correspondances de Langlands et de Jacquet-Langlands locales. Plus pr\'ecis\'ement, si $K$ est un corps $p$-adique de corps r\'esiduel $\FM_q,$ les composantes irr\'eductibles de la fibre sp\'eciale de l'espace sym\'etrique de Drinfeld pour $\GL_d(K)$ sont naturellement isomorphes \`a $\o X(w)$ et nous montrons dans \cite{RevMod} que les cycles proches du rev\^etement mod\'er\'e de Drinfeld sur une telle composante s'identifient naturellement \`a $Rj_*(\pi_*\L).$ Le th\'eor\`eme ci-dessus nous permet alors, toujours dans \cite{RevMod}, de calculer explicitement, et de mani\`ere purement locale, la cohomologie du rev\^etement mod\'er\'e de Drinfeld, et d'en d\'eduire en particulier qu'elle r\'ealise les correspondances de Langlands et Jacquet-Langlands pour les repr\'esentations elliptiques mod\'er\'ement ramifi\'ees de $\GL_d(K).$

\medskip

L'\'enonc\'e du th\'eor\`eme fait sens pour n'importe quel \'el\'ement de Coxeter $w$ d'un groupe r\'eductif $\Gb$ d\'efini sur $\FM_q,$ et nous conjecturons qu'il est vrai dans cette g\'en\'eralit\'e (voir \ref{R1S1}). La d\'emonstration du th\'eor\`eme se fait en trois \'etapes.

La premi\`ere \'etape de notre preuve utilise le calcul explicite de la normalisation de $\o X(w)$ dans $Y(w)$ d\^u \`a Bonnaf\'e et Rouquier dans \cite{bonnafe-rouquier-compact} pour se ramener au m\^eme \'enonc\'e pour la strate ``ouverte" $j_\Pb:X_\Pb(w):=\o X_\Pb(w)\ba \bigcup_{\Qb\subset\Pb}\o X_\Qb(w)\into \o X(w).$ Cette \'etape, d\'ecrite au paragraphe \ref{R1S2}, fonctionne pour tout $\Gb.$

La deuxi\`eme \'etape \'etudie la compactification partielle $X^\Pb(w):=X(w)\cup X_{\Pb}(w)\subset \o X(w)$ lorsque $\Pb$ est propre et maximal. On construit au paragraphe \ref{R1S3} un isomorphisme
$$
\Ub^F\ba X^\Pb(w)\simto X_\Lb(w_\Lb)\times\AM^1
$$
compatible avec l'isomorphisme de Lusztig $\Ub^F\ba X(w)\simto X_\Lb(w_\Lb)\times\GM_m$ de \cite{lusztig-coxeter}. Ici $\Lb$ d\'esigne le quotient r\'eductif de $\Pb$ est $w_\Lb$ un \'el\'ement de Coxeter de $\Lb.$ Pour cette construction, on utilise une description alternative de $\o X(w)$ comme \'eclat\'e d'espace projectif, valable seulement pour $\GL_d$ d\'eploy\'e. Dans le paragraphe \ref{ff}, on d\'eduit le th\'eor\`eme dans ce cas $\Pb$ maximal par un argument inspir\'e de Dudas \cite{dudas} et qui devrait fonctionner en toute g\'en\'eralit\'e, si on a un isomorphisme comme ci-dessus.

Enfin, la troisi\`eme \'etape est une r\'ecurrence sur le corang de $\Pb.$ Pour $\Gb$ quelconque, on explique cette r\'ecurrence au paragraphe \ref{R1S5}, sous l'hypoth\`ese que le cas de corang 1 est connu pour tout sous-groupe de Levi de $G.$

 \medskip

{\em Organisation de l'article.} Dans la section 2, nous rappelons les pr\'eliminaires sur les vari\'et\'es de Deligne-Lusztig, notamment la construction de Bonnaf\'e et Rouquier \cite{bonnafe-rouquier-compact} et le morphisme de quotient de Lusztig \cite{lusztig-coxeter}. Ensuite, on d\'ecrit en toute g\'en\'eralit\'e la premi\`ere et la troisi\`eme \'etapes mentionn\'ees ci-dessus. Dans la section 4, nous \'etudions le cas associ\'e \`a $\GL_d.$ Finalement, nous d\'emontrons notre th\'eor\`eme dans la section 5.

\medskip
\textbf{Remerciements:} Je remercie profond\'ement mon directeur de th\`ese Jean-Fran\c cois Dat pour les nombreuses discussions et ses constants encouragements pendant ces ann\'ees. Je remercie Zhi Jiang pour les conversations sur l'\'eclatement. Enfin, je remercie le referee anonyme pour ses suggestions qui ont permis d'en am\'eliorer consid\'erablement la r\'edaction.

\section{G\'en\'eralit\'es}
Dans cette partie, on rappelle tout d'abord la d\'efinition des vari\'et\'es de Deligne-Lusztig \cite{deligne-lusztig} et la construction de Bonnaf\'e et Rouquier \cite{bonnafe-rouquier-compact}. Ensuite, dans le cas de Coxeter, on rappelle certains r\'esultats de Lusztig \cite{lusztig-coxeter}.

\subsection{Pr\'eliminaires}\label{new5}
\ali Nous fixons un groupe r\'eductif connexe $\Gb$ d\'efini sur une cl\^oture alg\'ebrique $\o{\FM}_q$ du corps fini $\FM_q$. Nous supposons de plus que $\Gb$ est obtenu par extension des scalaires de $\Gb_0$ sur $\FM_q$, et nous notons $F:\Gb\to\Gb$ l'endomorphisme de Frobenius correspondant.

Fixons un \borel $F$-stable $\Bb$ de $\Gb$, un tore maximal $F$-stable $\Tb$ de $\Bb$ et notons $\Ub$ le radical unipotent de $\Bb$. Notons $W=N_\Gb(\Tb)/ \Tb$ le groupe de Weyl de $\Gb$ relativement \`a $\Tb$, $X(\Tb)$ (resp. $Y(\Tb)$) le r\'eseau des caract\`eres (resp. des sous-groupes \`a un param\`etre) de $\Tb$, $\Phi$ (resp. $\Phi^\vee$) le syst\`eme de racines (resp. coracines) de $\Gb$ relativement \`a $\Tb$, $\Delta$ (resp. $\Delta^\vee$) la base de $\Phi$ (resp. $\Phi^\vee$) associ\'ee \`a $\Bb$ et $\Phi_+$ (resp. $\Phi_+^\vee$) l'ensemble des racines (resp. coracines) positives contenant $\Delta$ (resp. $\Delta^\vee$). En particulier, $\Phi$ est stable sous l'action de $F.$ Si $\a$ est une racine, $F(\a)$ est un multiple positif d'une unique racine que l'on notera $\phi(\a),$ d\'efinissant une bijection $\phi:\Phi\to\Phi$ qui stablise $\D$ et $\Phi_+.$

Pour $\alpha\in\Phi,$ on notera $\alpha^\vee$ sa coracine associ\'ee, $s_\alpha\in W$ la r\'eflexion par rapport \`a $\alpha$, $\Ub_\alpha$ le sous-groupe unipotent \`a un param\`etre normalis\'e par $\Tb$ associ\'e \`a $\alpha,$ et $\Gb_\a$ le sous-groupe de $\Gb$ engendr\'e par $\Ub_\a$ et $\Ub_{-\a}.$ Posons $S:=\{s_\alpha\;|\;\alpha\in\Delta\}.$ On d\'esignera $l:W\to\NM$ la fonction longueur relativement \`a $S$. D'apr\`es \cite[9.3]{springer}, il existe un ensemble de repr\'esentants $\{\dot{w}\}$ de $W$ dans $N_\Gb(\Tb)$ v\'erifiant la propri\'et\'e suivante: si $w=w_1w_2$ est tel que $l(w)=l(w_1)+l(w_2)$, alors $\dot{w}=\dot{w_1}\dot{w_2}.$

Pour un \'el\'ement $w$ du groupe de Weyl $W$, Deligne et Lusztig ont construit dans \cite{deligne-lusztig} deux vari\'et\'es $X(w)$ et $Y(\dot{w})$ (ou $X_\Gb(w)$ et $Y_{\Gb}(\wdo)$ s'il y a confusion possible) sur $\o{\FM}_q$ ainsi qu'un morphisme fini \'etale $Y(\dot{w})\to X(w)$ faisant de $X(w)$ un quotient de $Y(\dot{w})$ par l'action du groupe fini $\Tb^{wF}:=\{t\in\Tb \;| \; wF(t)w^{-1}=t\}$ (dans {\em loc. cit.}, la vari\'et\'e $Y(\dot{w})$ est not\'ee $\wt{X}(\dot{w})$). Rappelons ci-dessous leurs d\'efinitions:
\begin{align*}
Y(\dot{w})&:=Y_\Gb(\dot{w})=\{g\Ub\in\Gb/\Ub\;|\;g^{-1}F(g)\in\Ub\dot{w}\Ub\} \\
X(w)&:=X_\Gb(w)=\{g\Bb\in\Gb/\Bb\;|\;g^{-1}F(g)\in\Bb w\Bb\}.
\end{align*}
Le groupe fini $\Gb^F$ agit par multiplication \`a gauche sur les vari\'et\'es quasi-projectives $X(w)$ et $Y(\dot{w}).$ De plus, le groupe commutatif $\Tb^{wF}$ agit librement sur $Y(\dot{w})$ par multiplication \`a droite. Le morphisme $\pi_w:Y(\dot{w})\to X(w)$ induit par restriction \`a $Y(\dot{w})$ de la projection canonique $\Gb/\Ub\onto\Gb/\Bb$ s'identifie \`a quotienter par $\Tb^{wF},$ induisant ainsi un isomorphisme $\Gb^F$-\'equivariant $Y(\dot{w})/\Tb^{wF}\simto X(w).$ Les vari\'et\'es $X(w)$ et $Y(\dot{w})$ ainsi obtenues sont quasi-affines, lisses et purement de dimension $l(w).$

Dans leur travail, Deligne et Lusztig ont construit \'egalement une compactification de Bott-Samelson-Demazure-Hansen des vari\'et\'es $X(w)$ (\cite[9.10]{deligne-lusztig}) que nous rappelons ci-dessous.

Soit $w=s_{\alpha_1}\cdots s_{\alpha_r},$ une expression minimale de $w,$ par rapport \`a la fonction longueur. On notera $w_i:=s_{\a_1}\cdots s_{\a_{i-1}}s_{\a_{i+1}}\cdots s_{\a_{r-1}}\in W.$ Posons alors, suivant \cite{deligne-lusztig}
\begin{align*}
\o{X}(w):=\biggl\{(g_1\Bb,\ldots,g_{r+1}\Bb)&\in(\Gb/\Bb)^{r+1}\;|\;
g_{r+1}\Bb=F(g_1)\Bb,\\
&\quad\quad g_i^{-1}g_{i+1}\in\Bb s_{\alpha_i}\Bb\cup\Bb,\;\forall 1\leq i\leq r
\biggl\}.
\end{align*}

D'apr\`es \cite[1.2]{deligne-lusztig}, la vari\'et\'e $X(w)$ s'identifie \`a la sous vari\'et\'e ouverte
\[
\{(g_1\Bb,\ldots,g_{r+1}\Bb)\in(\Gb/\Bb)^{r+1}\;|\;g_{r+1}\Bb=F(g_1)\Bb,~ g_i^{-1}g_{i+1}\in\Bb s_{\alpha_i}\Bb,\;\forall 1\leq i\leq r\}
\]
de $\o{X}(w).$ De plus, $\o{X}(w)$ est une vari\'et\'e lisse projective et $\o{X}(w)\backslash X(w)=\bigcup_{1\leq i\leq r}\o{X(w_i)}$ est un diviseur \`a croisements normaux (\cite[Lemme 9.11]{deligne-lusztig}). Notons que
\ini\begin{equation}
\o{X}(w)=\mathop{\coprod_{x=x_1\cdots x_r\in W}}_{x_i\in\{1,s_{\alpha_i}\}}X(x).
\end{equation}

\ali Dans \cite{bonnafe-rouquier-compact}, Bonnaf\'e et Rouquier ont donn\'e une construction explicite de la normalisation de $\o{X}(w)$ dans $Y(\dot{w}),$ not\'e $\o{Y}(\wdo).$ C'est l'unique vari\'et\'e normale $Z$ contenant $Y(\wdo)$ comme sous-vari\'et\'e ouverte dense et munie d'un morphisme fini $\o{\pi}_w:Z\onto \o{X}(w)$ prolongeant $\pi_w.$ Rappelons ci-dessous leurs constructions.

Tout d'abord, on peut supposer que le groupe d\'eriv\'e de $\Gb$ est simplement connexe ({\em cf.} \cite{bonnafe-rouquier-compact}). Ceci implique que $\Gb_\a\simto \SL_2$ et que $\a^\vee$ est injective pour toute racine $\a$. Si $1\leq i\leq r,$ il existe un unique $\l_i\in Y(\Tb)$ et un unique $m_i\in\ZM$ v\'erifiant les trois propri\'et\'es suivantes:
$$
\begin{cases}
\l_i-wF(\l_i)=m_i~s_1\cdots s_{i-1}(\a_i^\vee),\\
m_i>0,\\
Y(\Tb)/\ZM \l_i \text{ est sans torsion.}
\end{cases}
$$
Bonnaf\'e et Rouquier d\'efinissent une fonction $\varphi_\a:\Gb_\a\Ub\to\AM^1$ satisfaisant les propri\'et\'es de {\em loc. cit.} Prop. 2.2. Comme eux, posons
$$
\wt{\UC}(w):=\{\gb:=(g_1\Ub,\ldots,g_{r+1}\Ub)\in (\Gb/\Ub)^{r+1}~|~ \forall 1\leq i\leq r,~g_i^{-1}g_{i+1}\in\Gb_{\a_i}\Ub\},
$$
et notons
\begin{align*}
\wt{v}_w:\wt{\UC}(w) \To{} & \Gb/\Ub\times\Gb/\Ub\\
(g_1\Ub,\ldots,g_{r+1}\Ub)\mapsto &  (g_1\Ub,g_{r+1}\Ub).
\end{align*}

Soit $\mb:=(m_1,\ldots,m_r)$ une suite de $r$ entiers, notons

$$
\wt\UC_\mb(w):=\{(\gb,\xib)\in\wt\UC(w)\times\AM^r~|~\varphi_{\a_i}(g^{-1}_ig_{i+1})=\x_i^{m_i}\}
$$
et de m\^eme
$$
\wt\UC_\mb^I(w):=\{(\gb,\xib)\in \wt\UC_\mb(w)~|~\forall i\in I,~\x_i=0\}
$$
pour un sous-ensemble $I$ de $\{1,\ldots,r\}.$ Notons $$\G_F:=\{(g_1\Ub,g_2\Ub)\in \Gb/\Ub\times\Gb/\Ub~|~g_2\Ub=F(g_1)\Ub\}$$ le graphe du morphisme de Frobenius. Consid\'erons la vari\'et\'e
$$
\wt Y(w):=\{(\gb,\xib)\in \wt\UC_\mb(w)~|~ \wt v_w(\gb)\in \G_F\}
$$
et ses sous-vari\'et\'es localement ferm\'ees param\'etr\'ees par les sous-ensembles de $\{1,\ldots,r\}$:
$$
\wt Y_I(w):=\{(\gb,\xib)\in \wt\UC_\mb^I(w)~|~ \wt v_w(\gb)\in \G_F\}.
$$
Bonnaf\'e et Rouquier d\'efinissent une action naturelle de $\Tb^{wF}\times(\GM_m)^r$ sur la vari\'et\'e $\wt Y(w),$ faisant $\o X(w)$ le quotient de $\wt Y(w)$ par cette action, {\em cf.} \cite[Prop. 2.6]{bonnafe-rouquier-compact}, et $\wt Y_I(w)$ est stable sous cette action. Ils d\'emontrent le th\'eor\`eme suivant:

\begin{theo} \textsl{(\cite[Th\'eor\`eme. 1.2 (b)]{bonnafe-rouquier-compact})} La vari\'et\'e $\o Y(\wdo)=\wt Y(w)/(\GM_m)^r$ est normale et elle est munie d'une action de $\Tb^{wF}$ prolongeant l'action sur $Y(\wdo)$ telle que le morphisme de quotient $\o Y(\wdo)=\wt Y(w)/(\GM_m)^r\onto \o X(w)$ induit un isomorphisme $\o Y(\wdo)/\Tb^{wF}\simto \o X(w);$ autrement dit, elle est la normalisation de $\o{X}(w)$ dans $Y(\dot{w}).$
\end{theo}

Pour la commodit\'e du lecteur, on donne quelques propri\'et\'es faciles de la vari\'et\'e $\wt Y(w)$ qui seront utilis\'ees dans la preuve du lemme \ref{red 3}.
\begin{lemme}
Pour $I$ un sous-ensemble de $\{1,\ldots,r\},$ la vari\'et\'e $\wt \UC_\mb^I(w)$ est lisse, de dimension $2r+\dim \Gb/\Ub-|I|.$
\end{lemme}
\begin{preuve}
Notons $w(i):=s_{\a_1}\cdots s_{\a_{i}},$ $\mb_i:=(m_1,\ldots,m_i)$ et $I_i:=\{1,\ldots,i\}\cap I.$ On dispose d'une suite de morphismes canoniques ({\em cf.} \cite[Prop. 2.3]{bonnafe-rouquier-compact})
$$
\wt \UC_\mb^I(w)=\wt\UC_{\mb_r}^{I_r}(w(r))\to\wt\UC_{\mb_{r-1}}^{I_{r-1}}(w(r-1)) \to\cdots\to \wt\UC_{\mb_1}^{I_1}(w(1))\to \Gb/\Ub
$$
consistant \`a chaque \'etape \`a oublier le dernier terme de $\gb$ et $\xib.$ Lorsque $i\not\in I,$ la fibre du morphisme $\wt\UC_{\mb_i}^{I_i}(w(i))\to\wt\UC_{\mb_{i-1}}^{I_{i-1}}(w(i-1))$ est isomorphe \`a $\UC_{\a_i,m_i}=\{(g,\x)\in\Gb_{\a_i}\Ub/\Ub\times\AM^1~|~\varphi_{\a_i}(g)=\x^{m_i}\}.$ D'apr\`es \cite[Prop. 2.3]{bonnafe-rouquier-compact},
$$
\UC_{\a_i,m_i}\simto\{(x,y,\x)\in\AM^3~|~(x,y)\neq(0,0)\text{ et } y=\x^{m_i}\}\simto \AM^2\ba\{(0,0)\}.
$$
Lorsque $i\in I,$ la fibre du morphisme $\wt\UC_{\mb_i}^{I_i}(w(i))\to\wt\UC_{\mb_{i-1}}^{I_{i-1}}(w(i-1))$ est isomorphe \`a $\{(g,\x)\in\Gb_{\a_i}\Ub/\Ub\times\AM^1~|~\varphi_{\a_i}(g)=0\}\simto \AM^1\ba \{0\}.$ Donc c'est une suite des fibrations successives de fibres successivement isomorphes \`a des vari\'et\'es lisses (de dimension $1$ ou $2$), d'o\`u l'\'enonc\'e du lemme.

\end{preuve}

\begin{coro}
La vari\'et\'e $\wt Y_I(w)$ est lisse, purement de dimension $2r-|I|.$
\end{coro}
\begin{preuve}
Ceci d\'ecoule directement du lemme pr\'ec\'edent et \cite[Lemme 2.5]{bonnafe-rouquier-compact}.
\end{preuve}

\begin{coro}\label{dcn}
La vari\'et\'e $\wt Y_\vide(w)$ est une sous-vari\'et\'e ouverte dense dans $\wt Y(w)$ dont le compl\'ementaire est un diviseur \`a croisements normaux.
\end{coro}
\begin{preuve}
La premi\`ere assertion est dans \cite[Page 634-635]{bonnafe-rouquier-compact}. La deuxi\`eme assertion repose sur le corollaire pr\'ec\'edent. En effet, pour $1\leq i\leq r,$ notons $\wt Y_i(w):=\wt Y_{\{i\}}(w).$ Chaque $\wt Y_i(w)$ est un diviseur de $\wt{Y}(w)$, et
$$
\wt Y(w)\ba \wt Y_\vide(w)=\bigcup_{1\leq i\leq r}\wt Y_i(w).
$$
Gr\^ace au corollaire pr\'ec\'edent, $$\wt Y_I(w)=\bigcap_{i\in I}\wt Y_i(w)$$ est de codimension $|I|$ dans $\wt Y(w),$ $\forall I\subset\{1,\ldots,r\}.$
\end{preuve}

\subsection{Les orbites de Coxeter}\label{new9}
Soient $n=|\D/\phi|$ et $[\D/\phi]:=\{\alpha_1,\ldots,\alpha_n\}$ un syst\`eme de repr\'esentants de $\D/\phi.$ Notons $c=s_{\a_1}\cdots s_{\a_n}$ (ou $c_{\Gb}$ s'il y a confusion possible) un \'el\'ement de Coxeter, $w_\Delta$ l'\'el\'ement de $W$ de longueur maximale. Soit $I$ un sous-ensemble de racines simples stable sous $\phi,$ nous notons $W_I$ le sous-groupe de $W$ engendr\'e par les $\{s_\alpha\}_{\alpha\in I},$ $\Pb_I$ le \para $\Bb W_I\Bb$ de $\Gb.$ Posons $\Tb_I$ la composante connexe de l'\'el\'ement neutre de $\cap_{\alpha\in I}\Ker\alpha$, $\Lb_I:=Z_{\Gb}(\Tb_I)$ l'unique composante de Levi de $\Pb_I$ contenant $\Tb,$ $\Bb_I$ le \borel $\Bb\cap\Lb_I$ de $\Lb_I,$ $\Ub_I$ le radical unipotent de $\Pb_I,$ et $\Vb_I$ le radical unipotent de $\Bb_I.$ Le groupe de Weyl $W_{\Lb_I}$ de $\Lb_I$ associ\'e \`a $\Tb$ s'identifie \`a $W_I.$ On notera les groupes de points $F$-stables par les caract\`eres non \'epaissis correspondants $G, ~ P_I,~ L_I,~ U_I,~V_I,~\ldots$

Pour un sous-ensemble propre $I$ de $\D$ stable sous $\phi,$ on peut consid\'erer la vari\'et\'e de Deligne-Lusztig $X_{\Lb_I}(c_I)$ associ\'ee au groupe r\'eductif $\Lb_I$ et l'\'el\'ement de Coxeter $c_I$ de $W_{\Lb_I}$ obtenu \`a partir de $c$ en ne gardant que les r\'eflexions simples de $I.$ Notons $X_I:=X(x_1\cdots x_n)$, o\`u $x_i\in\{1,s_{\a_i}\}$ et $x_i=s_{\a_i}$ si et seulement si $\a_i\in I.$ D'apr\`es Lusztig \cite{lusztig-finiteness}, les vari\'et\'es $X_{\Lb_I}(c_I)$ et $X_I$ sont reli\'ees par la propri\'et\'e suivante:

\begin{prop}\textsl{(\cite[1.17]{lusztig-coxeter}, \cite[Lemme 3]{lusztig-finiteness}, voir aussi \cite[Prop. 3.3]{bonnafe-rouquier-coxeter})}\label{ee}
Sous l'hypoth\`ese comme plus haut, on a des isomorphismes canoniques:
\begin{align*}
G/U_I\times_{L_I}X_{\Lb_I}(c_I)&\simto X_I\\
(gU_I,h\Bb_I)&\mapsto gh\Bb.
\end{align*}
\end{prop}
\medskip

Comme $c_I$ est un \'el\'ement de Coxeter de $W_{\Lb_I},$ la vari\'et\'e $X_{\Lb_I}(c_I)$ est irr\'eductible d'apr\`es \cite[Prop. 4.8]{lusztig-coxeter}. Alors, $X_I$ est une union disjointe de ses composantes irr\'eductibles chacune isomorphe \`a $X_{\Lb_I}(c_I).$ Notons alors $C_I$ la composante irr\'eductible de $X_I$ fix\'ee par $U_I,$ d'apr\`es \cite[1.17]{lusztig-coxeter},
$$
C_I=\{(g_1\Bb,\ldots,g_{n+1}\Bb)\in X_I ~|~ g_1\Bb g^{-1}_1\subset\Pb_I\},
$$
munie d'un isomorphisme $L_I$-\'equivariant avec $X_{\Lb_I}(c_I).$ Notons ensuite $\o C_I$ son adh\'erence dans $\o X(c),$ alors
$$
\o C_I=\{(g_1\Bb,\ldots,g_{n+1}\Bb)\in \o X_I ~|~ g_1\Bb g^{-1}_1\subset\Pb_I\}.
$$
On en d\'eduit que $U_I$ agit trivialement sur $\o C_I,$ et $C_I$ n'a pas de points fixes sous un sous-groupe unipotent de $L_I,$ comme $C_I$ est une vari\'et\'e de Coxeter.

Plus g\'en\'eralement, on peut d\'efinir la stratification index\'ee par les \paras $F$-stables mentionn\'ee dans l'introduction. Pour $\Pb$ un \para $F$-stable conjugu\'e \`a $\Pb_I,$ on note $X_\Pb:=g\cdot C_I,$ o\`u $g\in G$ est tel que $\Pb=g\Pb_I g^{-1}.$ Alors, $X_I$ est une union disjointe de $X_\Pb$ o\`u $\Pb$ parcourt l'ensemble des \paras $F$-stables conjugu\'es \`a $\Pb_I.$ La vari\'et\'e $X_\Pb$ est localement ferm\'ee, irr\'eductible, et son adh\'erence est
$$
\o X_\Pb=\bigcup_{\Qb\subset\Pb}X_\Qb.
$$
Notons $U_\Pb$ le sous-groupe unipotent des points rationnels du radical unipotent de $\Pb,$ alors $U_\Pb$ agit trivialement sur $X_\Pb$ et $\o X_\Pb.$ On d\'emontre que l'ensemble des points fixes de $\o X(c)$ sous l'action de $U_\Pb$ s'identifie \`a $\o X_{\Pb}.$ En effet, si $U_\Pb$ a des points fixes dans une strate $X_\Qb,$ alors $U_\Pb$ normalise $\Qb,$ donc $U_\Pb$ est contenu dans $\Qb.$ Par ailleurs, comme $X_\Qb$ n'a pas de points fixes sous un sous-groupe unipotent de son quotient de Levi $L_\Qb,$ on a $U_\Pb$ contenu dans $U_\Qb,$ donc $\Pb$ contient $\Qb$ et $X_\Qb\subset \o X_\Pb.$ En particulier, on sait alors que $\o C_I=\o X(c)^{U_I}.$

\ali
Lorsque les vari\'et\'es de Deligne-Lusztig sont associ\'ees \`a des \'el\'ements de Coxeter, Lusztig a construit dans \cite{lusztig-coxeter} leurs quotients par $U$ et $U_I.$ Rappelons ci-dessous leurs constructions. Notons comme dans \cite{bonnafe-rouquier-coxeter},
\ini\begin{equation}\label{cond de lusztig}
X'(c)=\big\{u\in\Ub\;|\;u^{-1}F(u)\in(\Ub_{-w_\Delta(\alpha_1)}\backslash\{1\})\times\cdots(\Ub_{-w_\Delta(\alpha_n)}\backslash\{1\})\big\}
\end{equation}

Tout d'abord, Lusztig a d\'emontr\'e le th\'eor\`eme suivant:
\begin{theo}\textsl{(\cite[2.5, 2.6]{lusztig-coxeter})}\label{lusztig thm}
\begin{description}

\item [(a)] $X(c)\subset\Bb w_\Delta\cdot\Bb/\Bb$.

\item [(b)] Le morphisme
\begin{align*}
L:X'(c)&\longrightarrow X(c)\\
u&\longmapsto uw_\Delta\cdot\Bb
\end{align*}
est un isomorphisme de vari\'et\'es.
\end{description}
\end{theo}

\begin{rem}
En faisant agir sur $X'(c)$ le $p$-groupe fini $U$ par multiplication \`a gauche, et le groupe commutatif $\Tb^{F}$ par conjugaison, l'isomorphisme $L$ est $B$-\'equivariant.
\end{rem}


\begin{theo}\textsl{(\cite[Corollaries 2.7, 2.10]{lusztig-coxeter})}
\begin{description}
\item [(a)] L'isomorphisme $L$ dans $(b)$ du th\'eor\`eme pr\'ec\'edent induit un isomorphisme
\[
(\GM_m)^n=(\Ub_{-w_\Delta(\alpha_1)}\backslash\{1\})\times\cdots\times(\Ub_{-w_\Delta(\alpha_n)}\backslash\{1\})\simto U\ba X(c).
\]

\item [(b)] On a un morphisme naturel induit par l'isomorphisme dans (a):
$$
U_I\ba X(c)\onto U\ba X(c)\simto(\Ub_{-w_\Delta(\alpha_1)}\backslash\{1\})\times\cdots\times(\Ub_{-w_\Delta(\alpha_n)}\backslash\{1\})\onto \Ub_{-w_\Delta(\alpha_i)}=\GM_m.
$$
Ce morphisme induit un isomorphisme $V_I$-\'equivariant
\[
U_I\ba X(c)\simto X_{\Lb_I}(c_I)\times\GM_m,
\]
ainsi qu'un isomorphisme $L_I$-\'equivariant de cohomologies:
$$
R\G(U_I\ba X(c),\L)\simto R\G(X_{\Lb_I}(c_I)\times\GM_m,\L).
$$
\end{description}
\end{theo}

\section{La compactification partielle}

\subsection{\'Enonc\'e de la conjecture A}\label{R1S1}
On utilise les notations de la section \ref{new9}, et on d\'esigne $i_I:\o C_I\into \o X(c)$ l'inclusion naturelle pour tout $I\subsetneqq \D$ stable sous $\phi.$ Nous avons le diagramme commutatif suivant:
$$
\xymatrix{
Y(\cdo)\ar[d]^{\pi} \ar@{^(->}[r]^{j'} &\o{Y}(\cdo) \ar@{<-^)}[r]^{i'_I} \ar[d]_{\o{\pi}} & \o{\pi}^{-1}(\o{C}_I) \ar[d]_{\o{\pi}_I}\\
X(c)\ar@{^(->}[r]^j & \o X(c) \ar@{<-^)}[r]^{i_I} & \o{C}_I
}
$$
o\`u $\o{\pi}^{-1}(\o{C}_I):=(\o{Y}(\cdo)\times_{\o{X}(c)}\o{C}_I)_{\re}$ et $\o{\pi}_I=\o{\pi}|_{\o{\pi}^{-1}(\o{C}_I)}.$ On a formul\'e dans l'introduction la conjecture suivante:

\noindent{\bf Conjecture A.} {\it
Le morphisme de restriction:
$$
R\G(X(c),\pi_*\L)=R\G(\o X(c) ,Rj_*(\pi_*\L))\To{\res.}R\G(\o{C}_I,i^*_IRj_*(\pi_*\L))
$$
induit un isomorphisme
\begin{equation*}
R\G(X(c),\pi_*\L)^{U_I}\simto R\G(\o{C}_I,i^*_IRj_*(\pi_*\L)),
\end{equation*}
o\`u $\L=\ZM/\ell^m\ZM,$ pour un nombre premier $\ell\neq p.$}

\subsection{Conjecture A $\Leftrightarrow$ Conjecture A'}\label{R1S2}

Dans cette partie, on d\'emontre que la conjecture A \'equivaut \`a la conjecture suivante:

\noindent{\bf Conjecture A'.} {\it Pour tout $I\subsetneqq\D$ stable sous $\phi,$ le morphisme de restriction induit un isomorphisme:
$$
R\G(Y(\cdo),\L)^{U_I}\simto R\G(\o{\pi}^{-1}(C_I),j'^{*}_IRj'_*\L),
$$
o\`u $j'_I$ d\'esigne le compos\'e $\o\pi^{-1}(C_I)\into\o\pi^{-1}(\o C_I)\To{i'_I}\o Y(\cdo).$
}

L'\'equivalence entre les conjectures A et A' d\'ecoule directement des lemmes \ref{red 2} et \ref{red 3}. Remarquons que l'on s'est d\'ebarrass\'e de la compactification $\o{C}_I,$ et on se ram\`ene \`a \'etudier la strate ouverte $C_I.$

\begin{lemme}\label{red 2}
Le morphisme de changement de base induit un isomorphisme
$$
R\G(\o{C}_I,i^*_IRj_*(\pi_*\L))\To{\cong} R\G(\o{\pi}^{-1}(\o{C}_I),i'^{*}_I Rj'_*\L)
$$
compatible avec les morphismes de restriction.
\end{lemme}
\begin{preuve}
Notons que $\o{\pi}$ est un morphisme fini, il r\'esulte du th\'eor\`eme de changement de base pour les morphismes finis et l'invariance topologique du topos \'etale ({\em cf.} \cite[Exp. VIII Thm. 1.1]{SGA4-2}) que $i^*_IRj_*(\pi_*\L)=\o{\pi}_{I*}i'^{*}_I Rj'_*\L.$
\end{preuve}

\begin{lemme}\label{red 3}
Consid\'erons la normalisation $\o Y(\cdo)$ et le diagramme suivant:
$$\xymatrix{
Y(\cdo)\ar[r]^{j'} &\o Y(\cdo)&\o{\pi}^{-1}(\o{C}_I) \ar[l]_{i'_I}\\
& &  \o{\pi}^{-1}(C_I)\ar[u]\ar[ul]^{j'_I}}
$$
Le morphisme canonique
$$
i'_{I,*}i'^*_I Rj'_*\L\To{} Rj'_{I,*}j'^*_I Rj'_*\L
$$
est un isomorphisme.
\end{lemme}

\begin{preuve}
Rappelons la construction explicite de $\o Y(\cdo)$ dans \cite{bonnafe-rouquier-compact}.

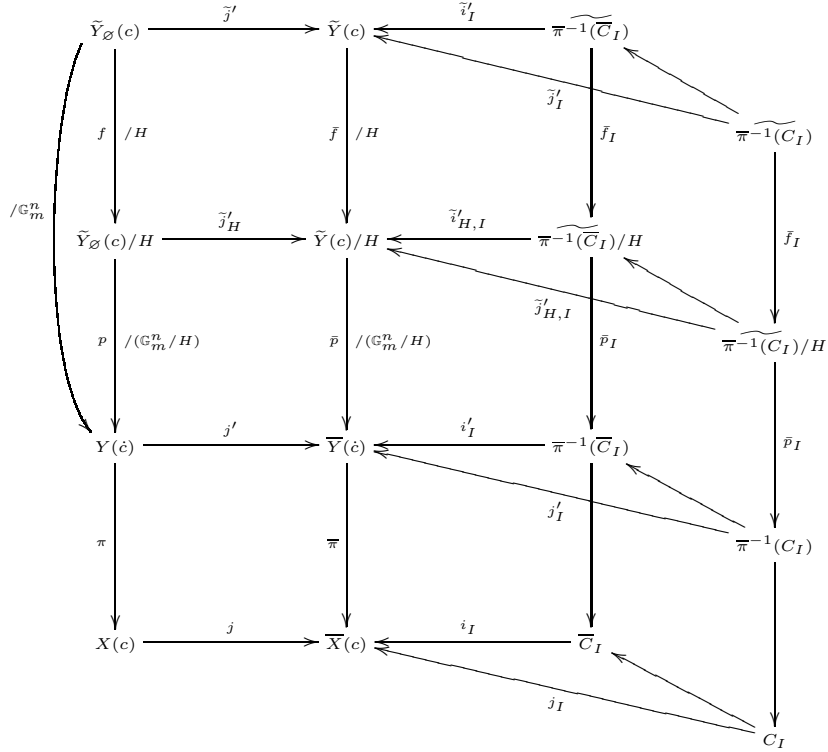
$$\tiny\xymatrix{
\wt Y_\vide(c)\ar@(lu,l)[dddd]_{/\GM_m^{n}}\ar[dd]_{f}^{/H}\ar[rr]^{\wt{j}'} & &\wt Y(c)\ar[dd]^{/H}_{\fba} & &\wt{\o{\pi}^{-1}(\o{C}_I)}\ar[ll]_{\wt{i}'_I}\ar[dd]^{\fba_I}&\\
& & & & &\wt{\o{\pi}^{-1}(C_I)}\ar[ulll]^{\wt{j}'_I}\ar[ul]\ar[dd]^{\fba_I}\\
\wt Y_\vide(c)/H\ar[dd]_{p}^{/(\GM_m^{n}/H)}\ar[rr]^{\wt{j}'_H} & &\wt Y(c)/H\ar[dd]_{\pba}^{/(\GM_m^{n}/H)}& &\wt{\o{\pi}^{-1}(\o{C}_I)}/H\ar[dd]^{\pba_I}\ar[ll]_{\wt{i}'_{H,I}} &\\
& & & & & \wt{\o{\pi}^{-1}(C_I)}/H\ar[ulll]^{\wt{j}'_{H,I}}\ar[ul]\ar[dd]^{\pba_I}\\
Y(\cdo) \ar[dd]_{\pi} \ar[rr]^{j'} & &\o{Y}(\cdo) \ar[dd]_{\o{\pi}}& &\o{\pi}^{-1}(\o{C}_I)\ar[dd]\ar[ll]_{i'_I}&\\
& & & & &\o{\pi}^{-1}(C_I)\ar[ulll]^{j'_I}\ar[ul]\ar[dd]\\
X(c)\ar[rr]^j & &\o{X}(c)& &\o{C}_I\ar[ll]_{i_I} &\\
& & & & & C_I\ar[ulll]^{j_I}\ar[ul]}
$$

Comme dans {\em loc. cit.} $\wt Y(c)$ est une vari\'et\'e lisse, purement de dimension $2n,$ munie d'une action de $\Tb^{cF}\times (\GM_m)^{n},$ faisant de $\o{X}(c)$ le quotient par cette action. L'application $\wt{j}':\wt Y_\vide(c)\into \wt Y(c)$ est une immersion ouverte dense, dont le compl\'ementaire est un diviseur \`a croisements normaux ({\em cf.} \ref{dcn}). Un argument similaire \`a l'\'etape 3 de la preuve du \cite[Thm. 2.2]{dat-lemma} nous fournit un isomorphisme:
\ini\begin{equation}\label{open 1}
\wt{i}'_{I,*}\wt{i}'^*_I R\wt{j}'_*\L\simto R\wt{j}'_{I,*}\wt{j}'^*_I R\wt{j}'_*\L.
\end{equation}
L'action du groupe $(\GM_m)^{n}$ sur $\wt Y(c)$ n'est pas libre en g\'en\'eral et il existe un sous-groupe fini $H$ (qui est not\'e $H_{\{1,\ldots,n\}}$ dans {\em loc. cit.}) de $\GM_m^{n}$ tel que $\GM_m^{n}/H$ agisse librement sur $\wt Y(c)/H$ (\cite[Prop. 2.7 (a)]{bonnafe-rouquier-compact}). Consid\'erons le diagramme commutatif suivant:
$$\xymatrix{
\wt Y_\vide(c) \ar[dd]_{f}^{/H}\ar[rr]^{\wt{j}'} & &\wt Y(c)\ar[dd]^{/H}_{\fba} & &\wt{\o{\pi}^{-1}(\o{C}_I)}\ar[ll]_{\wt{i}'_I}\ar[dd]^{\fba_I}&\\
& & & & &\wt{\o{\pi}^{-1}(C_I)}\ar[ulll]^{\wt{j}'_I}\ar[ul]\ar[dd]^{\fba_I}\\
\wt Y_\vide(c)/H\ar[rr]^{\wt{j}'_H} & &\wt Y(c)/H& &\wt{\o{\pi}^{-1}(\o{C}_I)}/H\ar[ll]_{\wt{i}'_{H,I}} &\\
& & & & & \wt{\o{\pi}^{-1}(C_I)}/H\ar[ulll]^{\wt{j}'_{H,I}}\ar[ul]}
$$

Le groupe fini $H$ agit librement sur $\wt Y_\vide(c)$ (\cite[Prop. 2.7 (b)]{bonnafe-rouquier-compact}), donc $f$ est un rev\^etement galoisien. Le faisceau \'etale $f_*\L$ sur $\wt Y_\vide(c)/H$ est un faisceau de $\L[H]$-modules, et le morphisme d'adjonction induit un isomorphisme $\L\simto (f_*\L)^H.$ Notons $R_H$ le foncteur d\'eriv\'e du foncteur des $H$-invariants. Comme $H$ est fini, prendre les $H$-invariants est une limite projective finie, donc commute avec les foncteurs pull-back. D'autre part, $R_H=R\Hom_{\L[H]}(\L,-)$ commute avec les images directes. On en d\'eduit que

$$\wt{i}'_{H,I,*}\wt{i}'^*_{H,I} R\wt{j}'_{H,*}\L=\wt{i}'_{H,I,*}\wt{i}'^*_{H,I} R\wt{j}'_{H,*}(f_*\L)^H=R_H(\fba_*\wt{i}'_{I,*}\wt{i}'^*_I R\wt{j}'_*\L)$$
et $$R\wt{j}'_{H,I,*}\wt{j}'^*_{H,I} R\wt{j}'_{H,*}\L=R_H(\fba_*R\wt{j}'_{I,*}\wt{j}'^*_I R\wt{j}'_*\L).$$
En plus, le morphisme $$\wt{i}'_{H,I,*}\wt{i}'^*_{H,I} R\wt{j}'_{H,*}\L\to R\wt{j}'_{H,I,*}\wt{j}'^*_{H,I} R\wt{j}'_{H,*}\L$$ s'identifie \`a $$R_H\fba_*(\wt{i}'_{I,*}\wt{i}'^*_I R\wt{j}'_*\L\to R\wt{j}'_{I,*}\wt{j}'^*_I R\wt{j}'_*\L).$$
D'apr\`es \ref{open 1}, on a donc
\ini\begin{equation}\label{open 2}
\wt{i}'_{H,I,*}\wt{i}'^*_{H,I} R\wt{j}'_{H,*}\L\simto R\wt{j}'_{H,I,*}\wt{j}'^*_{H,I} R\wt{j}'_{H,*}\L
\end{equation}

Consid\'erons ensuite le diagramme commutatif:
$$\xymatrix{
\wt Y_\vide(c)/H\ar[dd]_{p}^{/(\GM_m^{n}/H)}\ar[rr]^{\wt{j}'_H} & &\wt Y(c)/H\ar[dd]_{\pba}^{/(\GM_m^{n}/H)}& &\wt{\o{\pi}^{-1}(\o{C}_I)}/H\ar[dd]^{\pba_I}\ar[ll]_{\wt{i}'_{H,I}} &\\
& & & & & \wt{\o{\pi}^{-1}(C_I)}/H\ar[ulll]^{\wt{j}'_{H,I}}\ar[ul]\ar[dd]^{\pba_I}\\
Y(\cdo)  \ar[rr]^{j'} & &\o{Y}(\cdo) & &\o{\pi}^{-1}(\o{C}_I)\ar[ll]_{i'_I}&\\
& & & & &\o{\pi}^{-1}(C_I)\ar[ulll]^{j'_I}\ar[ul]\\
}$$
Le morphisme $\pba$ est le quotient par l'action libre de $(\GM_m)^{n}/H$, donc il est lisse. D'apr\`es le th\'eor\`eme de changement de base lisse et l'invariance topologique du topos \'etale (\cite[Exp. VIII]{SGA4-2}), on a $$\pba^*(i'_{I,*}i'^*_I Rj'_*\L)=\wt{i}'_{H,I,*}\wt{i}'^*_{H,I} R\wt{j}'_{H,*}\L$$
et
$$\pba^*(Rj'_{I,*}j'^*_I Rj'_*\L)=R\wt{j}'_{H,I,*}\wt{j}'^{*}_{H,I} R\wt{j}'_{H,*}\L.$$
En plus, le morphisme $$\wt{i}'_{H,I,*}\wt{i}'^*_{H,I} R\wt{j}'_{H,*}\L\To{} R\wt{j}'_{H,I,*}\wt{j}'^{*}_{H,I} R\wt{j}'_{H,*}\L$$ s'identifie \`a $$\pba^*(i'_{I,*}i'^*_I Rj'_*\L\To{} Rj'_{I,*}j'^*_I Rj'_*\L).$$
D'apr\`es \ref{open 2}, c'est un isomorphisme: $$\pba^*(i'_{I,*}i'^*_I Rj'_*\L\simto Rj'_{I,*}j'^*_I Rj'_*\L). $$ En vertu de la surjectivit\'e de $\pba$, on en d\'eduit que le morphisme suivant: $$i'_{I,*}i'^*_I Rj'_*\L\To{} Rj'_{I,*}j'^*_I Rj'_*\L$$ est un isomorphisme.

\end{preuve}

\begin{lemme}\label{R1L2}
Soient $\Gb=\Gb_1\times \Gb_2,$ et $\cdo_{\Gb}=(\cdo_{\Gb_1},\cdo_{\Gb_2}).$ Supposons que la conjecture A' soit v\'erifi\'ee pour $\Gb_1$ et $\Gb_2,$ alors elle l'est aussi pour $\Gb.$
\end{lemme}
\begin{preuve}
Comme on a un isomorphisme $Y_{\Gb}(\cdo_{\Gb})\cong Y_{\Gb_1}(\cdo_{\Gb_1})\times Y_{\Gb_2}(\cdo_{\Gb_2})$ et idem pour la composante dans la compactification partielle, l'\'enonc\'e d\'ecoule de la formule de K\"unneth.
\end{preuve}

\subsection{R\'eduction au cas de codimension $1$}\label{R1S5}

Dans cette partie, on d\'emontre la conjecture A' sous certaines hypoth\`eses sur le cas de codimension $1.$ Si $I\subset \D$ stable sous $\phi$ est tel que $|I/\phi|=|\D/\phi|-1,$ on dit que c'est {\em un cas de codimension $1.$} Dans ce cas, $\dim X_{\Lb_I}(c_I)=\dim X(c)-1.$ On fait l'hypoth\`ese suivante:

\medskip

\noindent{\bf Hypoth\`ese.} {\it Pour toute composante de Levi $F$-stable $\Lb$ d'un \para propre $F$-stable de $\Gb,$ la conjecture A' est vraie pour tous les cas de codimension $1$ de $\Lb.$}

\medskip

Sous cette hypoth\`ese, on d\'emontre la conjecture A'. Autrement dit, la preuve pour les strates de codimension plus grande que un se ram\`ene \`a celle pour les strates de codimension \'egale \`a un. La d\'emonstration se fait par r\'ecurrence sur le rang semi-simple $\rg_{ss}(\Gb)$ de $\Gb.$ Tout d'abord, le cas $\rg_{ss}(\Gb)=1$ d\'ecoule de l'hypoth\`ese sur le cas de codimension $1.$ Soit $\Gb$ un groupe r\'eductif connexe de rang semi-simple $n.$ Supposons que la conjecture A' soit vraie pour tous les sous-groupes de Levi rationnels $\Lb$ de $\Gb$ tels que $\rg_{ss}(\Lb)<n.$ Soit $J$ un sous-ensemble des racines simples $\D$ de $\Gb$ stable sous $\phi.$ Par notre hypoth\`ese sur le cas de codimension $1,$ on peut supposer que $J$ est tel que $\rg_{ss}(\Lb_J)\leq \rg_{ss}(\Gb)-2.$ Choisissons alors un sous-ensemble $I\subset \D$ stable sous $\phi,$ contenant $J,$ et $\rg_{ss}(\Lb_I)=\rg_{ss}(\Gb)-1.$ On a alors le diagramme commutatif suivant:
$$\xymatrix{
& & \o{\pi}^{-1}(C_I)\ar[ld]_{j'_I}\ar[rd]^{j'^\circ_I}\ar[dd] & & \\
Y(\cdo)\ar[r]^{j'}\ar[dd]_{\pi} & \o{Y}(\cdo)\ar[dd]_{\o{\pi}} & & \o{\pi}^{-1}(\o{C}_I)\ar[ll]_{i'_I}\ar[dd] & \o{\pi}^{-1}(\o{C}_J)\ar[l]_{i'_{J\subset I}}\ar[dd]\ar@(rd,d)[lll]_{i'_J} \\
& & C_I\ar[ld]\ar[rd] & &\\
X(c) \ar[r] & \o{X}(c) & & \o{C}_I\ar[ll] & \o{C}_J\ar[l]
}
$$

\begin{lemme}\label{R1L1}
Les inclusions naturelles
$$\xymatrix{
\o{\pi}^{-1}(C_I)\ar@{^(->}^{j'^\circ_I}[r] & \o{\pi}^{-1}(\o{C}_I)\ar@{<-^)}^{i'_{J\subset I}}[r] & \o{\pi}^{-1}(\o{C}_J).
}
$$
induisent un isomorphisme:
$$
R\G(\o{\pi}^{-1}(C_I),\L)^{U_J/U_I}\simto R\G(\o{\pi}^{-1}(C_J),i'^*_{J\subset I}Rj'^\circ_{I *}\L|_{\o{\pi}^{-1}(C_J)}).
$$
\end{lemme}
\begin{preuve}
Notons $c_I$ l'\'el\'ement de Coxeter de $\Lb_I$ d\'efini dans \ref{new9}. D'apr\`es \ref{ee}, on a $C_I\cong X_{\Lb_I}(c_I).$ Consid\'erons le diagramme suivant:

$$\xymatrix{
& Y_{\Lb_I}(\cdo_I) \ar[ldd]_{i_{c_I}}\ar[d]^{\pi_{\Lb_I}}\ar[r]^{j'_{\Lb_I}}   &\o{Y}_{\Lb_I}(\cdo_I) \ar[d]^{\o{\pi}_{\Lb_I}} & \o{\pi}_{\Lb_I}^{-1}(C_J) \ar[l] _{i'_{\Lb_I,J}} \ar[d] \\
& C_I\cong X_{\Lb_I}(c_I)  \ar[r]^{j_{\Lb_I}}  & \o{X}_{\Lb_I}(c_I) & C_J\cong X_{\Lb_J}(c_J) \ar[l]_{i_{\Lb_I,J}}\\
\o{\pi}^{-1}(C_I) \ar[r] ^{j'^\circ_I} \ar[ru]^{\o{\pi}} & \o{\pi}^{-1}(\o{C}_I) \ar[ru]^{\o{\pi}}  & \o{\pi}^{-1}(C_J) \ar[ru] \ar[l] _{i'_{J\subset I}}
}$$
D'apr\`es la construction de Bonnaf\'e et Rouquier (\cite[ Thm. 1.2 (d)]{bonnafe-rouquier-compact}), il existe un morphisme canonique $i_{c_I}:Y_{\Lb_I}(\cdo_I)\onto \o{\pi}^{-1}(C_I)$ tel que $\o{\pi}^{-1}(C_I)$ soit le quotient de $Y_{\Lb_I}(\cdo_I)$ par un sous-groupe $N$ (qui est not\'e $N_{c_I}(Y_{c,c_I})$ dans {\em loc. cit.}). Par hypoth\`ese de r\'ecurrence sur $\Lb_I,$ on a un isomorphisme donn\'e par le morphisme de restriction:
$$
R\G(Y_{\Lb_I}(\cdo_I),\L)^{U_J/U_I}\simto R\G(\o{\pi}^{-1}_{\Lb_I}(C_J),i'^*_{\Lb_I,J}Rj'_{\Lb_I,*}\L).
$$
C'est \'equivalent \`a dire que le morphisme suivant induit par restriction est un isomorphisme
\ini\begin{equation}\label{R1E1}
R\G(C_I,\pi_{\Lb_I*}\L)^{U_J/U_I}\simto R\G(C_J,i^*_{\Lb_I,J}Rj_{\Lb_I,*}\pi_{\Lb_I *}\L).
\end{equation}

D'autre part, on a un diagramme commutatif
$$
\xymatrix{
R\G(\o{\pi}^{-1}(C_I),\L)^{U_J/U_I} \ar[r] \ar[d]^\cong &  R\G(\o{\pi}^{-1}(C_J),i'^*_{J\subset I}Rj'^\circ_{I *}\L|_{\o{\pi}^{-1}(C_J)})\ar[d]^\cong\\
R\G(C_I,\o{\pi}_*(\L_{\o{\pi}^{-1}(C_I)}))^{U_J/U_I} \ar[r]  &  R\G(C_J,i_{\Lb_I,J}^*Rj_{\Lb_I,*}\o{\pi}_*\L)}
$$
Notons que $\L_{\o{\pi}^{-1}(C_I)}=(i_{c_I,*}\L)^N$ et l'action de $U_J/U_I$ commute avec celle de $N,$ nous avons donc
\begin{align*}
R\G(C_I,\o{\pi}_*(\L_{\o{\pi}^{-1}(C_I)}))^{U_J/U_I}
&\simto R\G(C_I,\o{\pi}_*(i_{c_I,*}\L)^N)^{U_J/U_I}\\
&\simto R_N(R\G(C_I,\pi_{\Lb_I*}\L)^{U_J/U_I}\\
&\simto R_N(R\G(C_I,\pi_{\Lb_I*}\L)^{U_J/U_I}),
\end{align*}
et
\begin{align*}
R\G(C_J,i_{\Lb_I,J}^*Rj_{\Lb_I,*}\o{\pi}_*\L)
&\simto R\G(C_J,i_{\Lb_I,J}^*Rj_{\Lb_I,*}(\o{\pi}_*(i_{c_I,*}\L)^N))\\
&\simto R_N(R\G(C_J,i_{\Lb_I,J}^*Rj_{\Lb_I,*}\pi_{\Lb_I*}\L)).
\end{align*}
Notons que les isomorphismes ci-dessus sont compatibles avec les morphismes de restriction. En appliquant le foncteur d\'eriv\'e $R_N$ \`a l'isomorphisme \ref{R1E1}, on d\'eduit un isomorphisme
$$
R\G(C_I,\o{\pi}_*(\L_{\o{\pi}^{-1}(C_I)}))^{U_J/U_I}\simto R\G(C_J,i_{\Lb_I,J}^*Rj_{\Lb_I,*}\o{\pi}_*\L).
$$
Donc un isomorphisme
$$
R\G(\o{\pi}^{-1}(C_I),\L)^{U_J/U_I}\simto  R\G(\o{\pi}^{-1}(C_J),i'^*_{J\subset I}Rj'^\circ_{I *}\L|_{\o{\pi}^{-1}(C_J)}).
$$

\end{preuve}

Revenons \`a la preuve de la conjecture A'. Consid\'erons la compactification partielle $Y(\cdo)^I=Y(\cdo)\coprod\o{\pi}^{-1}(C_I).$ D'apr\`es \cite[Thm. 1.2 (a)]{bonnafe-rouquier-compact}, $(Y(\cdo)^I,\o{\pi}^{-1}(C_I))$ est un couple lisse de codimension $1.$ Le th\'eor\`eme de puret\'e relative (\cite[Exp. XVI]{SGA4-3}) nous fournit un triangle distingu\'e dans $D^b_c(\o{\pi}^{-1}(C_I),\L)$:
$$
\L\To{} j'^*_I Rj'_*\L\To{} \L(-1)[-1]\To{+1}.
$$
ainsi qu'un diagramme commutatif:
$$\xymatrix @C=4mm{
R\G\!(\!\o{\pi}^{\!-\!1}(C_I),\L)\ar[r]\ar[d]_{\res.} & R\G(\o{\pi}^{\!-\!1}(C_I),j'^*_IRj'_*\L) \ar[r] \ar[d]_{\res.} & R\G(\o{\pi}^{\!-\!1}(C_I),\L(\!-\!1))[\!-\!1]\ar[d]_{\res.}\ar[r]& .\\
R\G(\o{\pi}^{\!-\!1}(C_J),i'^*_{J\!\subset\! I}Rj'^\circ_{I*}\L) \ar[r] & R\G(\o{\pi}^{\!-\!1}(C_J),i'^*_{J\! \subset\! I}Rj'^\circ_{I*}j'^*_IRj'_* \L) \ar[r] &
R\G(\o{\pi}^{\!-\!1}(C_J),i'^*_{J\!\subset\! I}Rj'^\circ_{I*}\L(\!-\!1))[\!-\!1]\ar[r] &.}$$
En vertu du lemme pr\'ec\'edent, ceci induit un isomorphisme:
$$
R\G(\o{\pi}^{-1}(C_I),j'^*_IRj'_*\L)^{U_J/U_I}\simto R\G(\o{\pi}^{-1}(C_J),i'^*_{J\subset I}Rj'^\circ_{I*}(j'^*_IRj'_*\L|_{\o{\pi}^{-1}(C_I)}))
$$
D'apr\`es le lemme \ref{red 3}, $i'^*_JRj'_*\L=i'^*_{J\subset I}i'^*_IRj'_*\L=i'^*_{J\subset I}Rj'^\circ_{I*}(j'^*_IRj'_*\L|_{\o{\pi}^{-1}(C_I)}).$ On en d\'eduit que

\begin{align*}
R\G(Y(\cdo),\L)^{U_J}&=(R\G(Y(\cdo),\L)^{U_I})^{U_J/U_I}\\
&\simto R\G(\o{\pi}^{-1}(C_I),j'^*_IRj'_*\L)^{U_J/U_I}\\
&\simto R\G(\o{\pi}^{-1}(C_J),i'^*_JRj'_*\L).
\end{align*}
Ceci ach\`eve le raisonnement par r\'ecurrence.


\section{Le cas o\`u $\Gb=\GL_{d}(\o{\FM}_q)$}\label{new8}
Dans cette section, nous d\'etaillons les constructions de la section 2 lorsque $\Gb=\GL_d(\o{\FM}_q)$ le groupe lin\'eaire de dimension $d.$
\subsection{L'espace de Drinfeld sur un corps fini}
\ali Soient $\Gb=\GL_d(\o{\FM}_q)$ et $F$ l'endomorphisme de Frobenius standard $(a_{i,j})\mapsto(a_{i,j}^q).$ Via cette isog\'enie, $\Gb$ admet une $\FM_q$-structure d\'eploy\'ee de sorte que $\Gb^F,$ le groupe des points fix\'es par $F,$ s'identifie \`a $\GL_d(\FM_q).$ Notons de plus que $\Bb$ le \borel de $\GL_d(\o{\FM}_q)$ des matrices triangulaires sup\'erieures, $\Tb$ le tore maximal des matrices diagonales. Le syst\`eme de racines $\Phi$ est de type $A_{d-1}$, et $\Delta=\{\alpha_1,\ldots,\alpha_{d-1}\}$ la base de $\Phi$ associ\'ee \`a $\Bb$. On num\'erotera $\Delta$ de telle sorte que pour $i\in\{1,\ldots,d-1\}$, le parabolique $\Pb_{\Delta\backslash\{\alpha_i\}}$ soit le stabilisateur d'un sous-espace de dimension $i$ de $\oFq^d.$ Le groupe de Weyl $W$ s'identifie alors au groupe sym\'etrique $\SG_d,$ et nous choisissons l'\'el\'ement de Coxeter $c=(1,\ldots,d)\in \SG_d.$ D\'esormais, on notera $T_d:=\Tb^{cF},$ et l'application $(a_{ij})_{1\leq i,j\leq d}\mapsto a_{11}$ nous fournit un isomorphisme $T_d\simto \FM_{q^d}^\times.$

La vari\'et\'e projective $\Gb/\Bb$ s'identifie \`a l'ensemble $\FC$ des drapeaux complets de l'espace vectoriel $\o{\FM}_q^d$:
$$
\FC=\big\{\{0\}=D_0\subset D_1\subset\cdots\subset D_d=\o{\FM}_q^d\;|\;\dim_{\o{\FM}_q} D_i=i\big\}.
$$
En effet, $\Gb$ agit transitivement sur $\FC$ et $\Bb$ s'identifie au stabilisateur du drapeau canonique
$$
\{0\}\subset \DG_1\subset \DG_2\subset\cdots\subset \DG_{d-1}\subset \o{\FM}^d_q
$$
o\`u $\DG_i$ est le sous-espace d\'efini par $X_i=\cdots =X_{d-1}=0, ~\forall 1\leq i\leq d-1,$ et $X_0,\ldots,X_{d-1}$ d\'esignent les coordonn\'ees sous base canonique de l'espace vectoriel $\oFq^d.$ Via cette description, un drapeau $D_\bullet$ appartient \`a $X(c)$ \ssi $D_i=D_1\oplus F(D_1)\oplus\cdots\oplus F^{i-1}(D_1)$ pour tout $i$ ({\em cf.} \cite[2.2]{deligne-lusztig}). On obtient ainsi que $X(c)$ s'identifie \`a la sous-vari\'et\'e $\Omega_{\FM_q}^{d-1}$ de
$\PM^{d-1}_{\FM_q}$ 
d\'efinie, dans les coordonn\'ees projectives $[X_0:\ldots:X_{d-1}]$,
par la non-nullit\'e du d\'eterminant $\det((X_i^{q^j})_{0\leq i,j\leq d-1}).$ Autrement dit, $\Omega_{\FM_q}^{d-1}$ peut s'exprimer comme le compl\'ementaire de tous les hyperplans
$\FM_q$-rationnels dans $\PM^{d-1}_{\FM_q}$. Nous fixons $\cdo$ l'\'el\'ement dont l'action sur la base canonique $\{e_i\}_{1\leq i\leq d}$ est donn\'ee par $\cdo(e_i)=e_{i+1},~\forall 1\leq i\leq d-1$ et $\cdo(e_{d})=e_1.$ La vari\'et\'e $Y(\cdo)$ est finie \'etale sur $\Omeq^{d-1}$ de groupe de Galois $T_d\cong\FM_{q^d}^\times.$ De plus, on peut identifier
$Y(\dot{c})$ avec la sous-vari\'et\'e $\DL^{d-1}$ de l'espace affine
$\AM^d_{\FM_q}=\Spec {\FM_q[X_0,\ldots,X_{d-1}]}$ 
d\'efinie par l'\'equation $\det((X_i^{q^j})_{0\leq i,j\leq
d-1})^{q-1}=(-1)^{d-1}$. En particulier, $\DL^{d-1}$ est une vari\'et\'e affine
et lisse.

\ali On d\'esigne $\o{\O}^{d-1}_{\FM_q}$ la compactification de $\Omega_{\FM_q}^{d-1}$ ({\em cf.} \ref{new5}). Il est bien connu que $\o{\O}^{d-1}_{\FM_q}$ peut se construire par une suite d'\'eclatements successifs de l'espace projectif $\PM^{d-1}_{\FM_q}$ ({\em cf.} \cite[4.1]{ito}). Rappelons que si $X$ est une vari\'et\'e et $\pi:\wt X\to X$ est un \'eclatement de $X$ le long d'une sous-vari\'et\'e ferm\'ee $V,$ pour une sous-vari\'et\'e ferm\'ee $W$ non contenue dans $V,$ le transform\'e strict $\wt W\subset \wt X$ de $W$ s'identifie \`a l'adh\'erence de Zariski de $\pi^{-1}(W\ba V)$ dans $\wt X.$

\begin{lem}\label{compact}
On peut obtenir $\o{\O}^{d-1}_{\FM_q}$ de la fa\c con suivante: \`a partir de l'espace projectif $\PM^{d-1}_{\FM_q},$ on \'eclate tous ses points $\FM_q$-rationnels, et puis on \'eclate le long des transform\'es stricts de toutes les droites $\FM_q$-rationnelles, et puis on \'eclate le long des transform\'es stricts de tous les plans $\FM_q$-rationnels et ainsi de suite.
\end{lem}
\begin{preuve}
Notons $Z_i\subset Y_0:=\PM^{d-1}_{\FM_q}$ la r\'eunion de toutes les sous-vari\'et\'es lin\'eaires $\FM_q$-rationnelles de dimension $i,~\forall i\in\{0,\ldots,d-2\}.$ Posons $Y_1$ l'\'eclatement de $Y_0$ le long de $Z_0,$ et $Z_i^{(1)}$ le transform\'e strict de $Z_i$ dans $Y_1$ pour $i\geq 1.$ Posons $Y_2$ l'\'eclatement de $Y_1$ le long de $Z^{(1)}_1,$ et $Z_i^{(2)}$ le transform\'e strict de $Z^{(1)}_i$ pour $i\geq 2.$ Construisons les vari\'et\'es projectives $Y_k$ et $Z_i^{(k)}$ par r\'ecurrence. Supposons que $Y_{k-1}$ et $Z_i^{(k-1)}$ soient construites, posons $p_k:Y_k\to Y_{k-1}$ l'\'eclatement de $Y_{k-1}$ le long de $Z_{k-1}^{(k-1)},$ et $Z_i^{(k)}$ le transform\'e strict de $Z_i^{(k-1)}$ dans $Y_k.$ Notons que pour $j\leq k,$ $Z_{j}^{(k)}$ est une union disjointe de ses composantes irr\'eductibles. Finalement, on obtient une vari\'et\'e projective $Y_{d-2}$ et un morphisme $f:=p_1\circ \cdots \circ p_{d-2}:Y_{d-2}\To{}Y_0.$ Autrement dit, on a une suite d'\'eclatements:
$$
Y_{d-2}\To{p_{d-2}} Y_{d-3} \To{p_{d-3}}\cdots\To{p_2}Y_1\To{p_1}Y_0.
$$

Notons $X(c):=\O^{d-1}_{\FM_q},$ et $\o{X}(c):=\o{\O}^{d-1}_{\FM_q}=X(c)\coprod(\bigcup_{1\leq k\leq d-1}\o{X(c_k)}),$ o\`u $c_k=s_{\a_1}\cdots s_{\a_{k-1}}s_{\a_{k+1}}\cdots s_{\a_{d-1}}.$ On rappelle que
\begin{align*}
\o{X}(c):=\biggl\{(g_1\Bb,\ldots,g_{d}\Bb)&\in(\Gb/\Bb)^{d}\;|\;
g_{d}\Bb=F(g_1)\Bb,\\
&\quad\quad g_i^{-1}g_{i+1}\in\Bb s_{\alpha_i}\Bb\cup\Bb,\;\forall 1\leq i\leq d-1
\biggl\}
\end{align*}
est une sous-vari\'et\'e projective de $(\Gb/\Bb)^d,$ et que $\o{X}(c)\backslash X(c)$ est une r\'eunion des sous-vari\'et\'es projectives
\begin{align*}
\o{X(c_k)}=&\biggl\{\gb\in \o{X}(c)~\big|~ g_i^{-1}g_{i+1}\in \Bb s_{\a_{i}}\Bb\cup\Bb ~\forall i\in\{1,\ldots,k-1\},\\
 &\quad\quad\quad\quad      g_k^{-1}g_{k+1}\in\Bb,~g^{-1}_{j}g_{j+1}\in\Bb s_{\a_j}\Bb\cup\Bb~ \forall j\in\{k+1,\ldots,d-1\}   \biggl\}\\
\end{align*}
Notons pour $i\geq1,$ $\DG_i=\{X_j=\cdots=X_{d-1}=0\}\subset \PM^{d-1},$ et d\'efinissons un morphisme
\begin{align*}
\phi_0:\o{X}(c)&\To{}Y_0\\
\gb:=(g_1\Bb,\ldots,g_d\Bb)&\longmapsto g_1\cdot \DG_1
\end{align*}
C'est un morphisme projectif birationnel $\PGL_d(\FM_q)$-\'equivariant induisant un isomorphisme sur $X(c)\simto \Omeq^{d-1}.$ L'image r\'eciproque
$$
\phi^{-1}_0(\DG_1)=\{\gb\in\o{X(c_1)}~|~g_1=(a_{i,j})_{1\leq i,j\leq d}\text{ avec $a_{11}\neq 0,$ et $a_{i,1}=0$ $\forall i\geq 2$}\}
$$
est une composante irr\'eductible de $\o{X(c_1)},$ donc un diviseur de $\o{X}(c).$ Les $\FM_q$-points rationnels de $\PM^{d-1}$ sont conjugu\'es sous l'action de $\PGL_d(\FM_q),$ et leurs images inverses sous $\phi_0$ constituent un diviseur de $\o{X}(c).$ D'apr\`es la propri\'et\'e universelle d'\'eclatement \cite[II 7.14]{hartshorne} et le fait que $\phi_0$ est $\PGL_d(\FM_q)$-\'equivariant, il existe un unique morphisme projectif $\PGL_d(\FM_q)$-\'equivariant $\phi_1:\o{X}(c)\To{} Y_1$ tel que $\phi_0=p_1\circ \phi_1$ et $\phi_0^{-1}(Z_0)=\o{X(c_1)}.$ Notons $\DG_2^\circ$ le compl\'ementaire dans $\DG_2$ de tous ses points rationnels, et $\DG_2^{\circ(1)}$ (resp. $\DG^{(1)}_2$) le transform\'e strict de $\DG^\circ_2$ (resp. $\DG_2$) dans $Y_1.$  Alors on a
\begin{align*}
\phi^{-1}_1(\DG_2^{\circ(1)})=\phi^{-1}_0(\DG^\circ_2)=&\{\gb\in\o{X(c_2)}~|~g_1=(a_{ij})_{1\leq i,j\leq d}\text{ avec $a_{11},a_{21}$ $\FM_q$-lin\'eairement}\\
 &\quad\quad\quad\quad    \text{ind\'ependants, et $a_{i1}=0$ $\forall i\in\{3,\ldots,d\}$}  \}.
\end{align*}
De plus, sous ces conditions, on a $a_{ij}=0$ lorsque $3\leq i\leq d, 1\leq j\leq 2.$ On en d\'eduit un isomorphisme ({\em cf.} \ref{q})
$$
\phi^{-1}_{1}(\DG_2^{\circ(1)})\simto \O^{1}_{\FM_q}\times \o{\O}^{d-3}_{\FM_q}.
$$
Ceci entra\^ine que $\phi^{-1}_{1}(\DG_2^{(1)})$ est une composante irr\'eductible de $\o{X(c_2)},$ donc $\phi^{-1}_1(Z_1^{(1)})=\o{X(c_2)}.$ De plus, on a $\phi^{-1}_1(Z_0^{(1)}\cap Z_1^{(1)})=\phi^{-1}_1(Z_0^{(1)})\cap\phi^{-1}_1(Z_1^{(1)})=\o{X(c_1)}\cap\o{X(c_2)}.$

Par r\'ecurrence, construisons un morphisme projectif birationnel $\phi_k:\o{X}(c)\To{} Y_k$ pour chaque $k\leq d-2.$ Supposons que $\phi_{k-1}$ soit construit. Notons $\DG_k^\circ$ le compl\'ementaire dans $\DG_k$ de la r\'eunion des tous ses sous-espaces lin\'eaires rationnels propres, et $\DG_k^{\circ(k-1)}$ (resp. $\DG_k^{(k-1)}$) le transform\'e strict de $\DG_k^\circ$ (resp. $\DG_k$) dans $Y_{k-1}.$ Le morphisme $\phi_{k-1}$ induit bien s\^ur un isomorphisme entre $\DG_k^\circ$ et $\DG_k^{\circ(k-1)}.$ Comme $\DG_k$ est de dimension $k-1,$ $\DG_k^{(k-1)}$ est une composante irr\'eductible de $Z^{(k-1)}_{k-1}.$ Alors on a
\begin{align*}
\phi^{-1}_{k-1}(\DG_k^{\circ(k-1)})=\phi^{-1}_0&(\DG_k^\circ)=\{\gb\in\o{X(c_k)}~|~g_1=(a_{ij})_{1\leq i,j\leq d}\text{ avec $a_{11},\ldots,a_{k1}$} \\
 &\text{$\FM_q$-lin\'eairement }  \text{ind\'ependants, et $a_{i1}=0$ $\forall i\in\{k+1,\ldots,d\}$}  \}.
\end{align*}
De plus, sous ces conditions, on a $a_{ij}=0$ lorsque $k+1\leq i\leq d, 1\leq j\leq k.$ On en d\'eduit un isomorphisme ({\em cf.} \ref{q})
$$
\phi^{-1}_{k-1}(\DG_k^{\circ(k-1)})\simto \O^{k-1}_{\FM_q}\times \o{\O}^{d-1-k}_{\FM_q}.
$$
Ceci entra\^ine que $\phi^{-1}_{k-1}(\DG_k^{(k-1)})$ est une composante irr\'eductible de $\o{X(c_k)},$ donc un diviseur irr\'eductible. Par la $\PGL_d(\FM_q)$-\'equivariance de $Z_{k-1}^{(k-1)},$ on sait que $\phi^{-1}_{k-1}(Z_{k-1}^{(k-1)})=\o{X(c_k)}.$ En vertu de la propri\'et\'e universelle d'\'eclatement, il existe un unique morphisme projectif birationnel $\PGL_d(\FM_q)$-\'equivariant $\phi_k:\o{X}(c)\To{} Y_k$ qui rend le diagramme suivant commutatif:
$$
\xymatrix{
\o{X}(c)\ar[rrd]^{\phi_{k-1}}\ar[rrddd]^{\phi_0} \ar[rr]^{\phi_k} & & Y_k  \ar[d]^{p_k} \\
& & Y_{k-1}\ar[d]^{p_{k-1}}\\
& & \vdots \ar[d]^{p_1}\\
& & Y_0
}
$$
De la m\^eme mani\`ere, on sait que $\phi^{-1}_k(Z_k^{(k)})=\o{X(c_{k+1})},$ et pour $0\leq i_1,\ldots,i_n\leq k$
$$
\phi^{-1}_k(Z^{(k)}_{i_1}\cap\cdots\cap Z^{(k)}_{i_n})=\phi^{-1}_k(Z^{(k)}_{i_1})\cap\cdots\cap\phi^{-1}_k(Z^{(k)}_{i_n})=\o{X(c_{i_1+1})}\cap\cdots\cap\o{X(c_{i_n+1})}.
$$

Finalement, on obtient un morphisme projectif birationnel $\phi_{d-2}:\o{X}(c)\to Y_{d-2}$ v\'erifiant que
\ini\begin{equation}\label{r}
\phi^{-1}_{d-2}(\bigcap_{i\in S} Z^{(d-2)}_{i})=\bigcap_{i\in S}\phi^{-1}_{d-2}(Z^{(d-2)}_{i})=\bigcap_{i\in S}\o{X(c_{i+1})},
\end{equation}
pour n'importe quel sous-ensemble $S$ de $\{0,\ldots,d-2\}.$

Montrons maintenant que $\phi_{d-2}$ est un isomorphisme. Notons tout d'abord que $\phi_{d-2}$ est un \'eclatement d'apr\`es \cite[II 7.17]{hartshorne}. Supposons que $\phi_{d-2}$ soit un \'eclatement de $Y_{d-2}$ le long d'un sous-sch\'ema ferm\'e irr\'eductible $Z\subset Y_{d-2}\ba f^{-1}( \Omeq^{d-1})=\bigcup Z_i^{(d-2)}$ de codimension $\geq 2,$ o\`u $f$ est le compos\'e d'\'eclatements $Y_{d-2}\to Y_0.$ On peut alors supposer que $Z$ est contenu dans certain $Z_i^{(d-2)}.$ D'apr\`es \ref{r}, $Z$ n'est pas contenu dans l'intersection de $Z_i^{(d-2)}$ avec les autres composantes $Z_j^{(d-2)}.$ D'autre part, par construction, on sait que $\phi_{d-2}|_{\o{X(c_{i+1})}}:\o{X(c_{i+1})}\to Z_i^{(d-2)}$ est un morphisme birationnel. Or, le nerf de la stratification de $\o{X}(c)\ba X(c)$ et celui de $Y_{d-2}\ba f^{-1}(\Omeq^{d-1})$ co\"incident avec l'immeuble de Tits associ\'e \`a $\GL_d(\FM_q).$ Il s'ensuit que un tel sous-sch\'ema $Z$ n'existe pas. Donc $\phi_{d-2}$ est un isomorphisme. Ceci termine la preuve du lemme.
\end{preuve}

\ali\label{q} Lorsque $I$ est un sous-ensemble de $\Delta$ tel que $\Delta\backslash I=\{\alpha_i\},$ $\Lb_I$ s'identifie au produit $\GL_i(\o{\FM}_q)\times\GL_{d-i}(\oFq),$ et $\Bb_I$ s'identifie au produit $\Bb_i\times\Bb_{d-i}$ des \borels standards de $\GL_i(\oFq)$ et $\GL_{d-i}(\oFq)$ respectivement. Le groupe de Weyl $W_{\Lb_I}$ isomorphe au produit $W_{\GL_i(\oFq)}\times W_{\GL_{d-i}(\oFq)}$, de sorte qu'il existe des \'el\'ements $\o{s}_{\alpha_k}\in W_{\GL_i(\oFq)}$ pour $1\leq k\leq i-1$ et $\o{s}_{\alpha_l}\in W_{\GL_{d-i}(\oFq)}$ pour $i+1\leq l\leq d-1$, tels que via cet isomorphisme $s_{\alpha_k}$ s'identifie \`a $(\o{s}_{\alpha_k},1)$ pour $1\leq k\leq i-1$ et $s_{\alpha_l}$ s'identifie \`a $(1,\o{s}_{\alpha_l})$ pour $i+1\leq l\leq d-1.$ Notons $c_1(i):=\o{s}_{\alpha_1}\cdots\o{s}_{\alpha_{i-1}}$ (resp. $c_2(d-i):=\o{s}_{\alpha_{i+1}}\cdots\o{s}_{\alpha_{d-1}}$) l'\'el\'ement de Coxeter $(1,\ldots,i)\in \SG_i=W_{\GL_i(\oFq)}$ (resp. $(1,\ldots,d-i)\in \SG_{d-i}=W_{\GL_{d-i}(\oFq)}$).

D\'esormais, on identifiera $X_{\Lb_I}(c_I)$ avec $\Omeq^{i-1}\times\Omeq^{d-1-i}$ de la fa\c con suivante:

\begin{align*}
X_{\Lb_I}(c_I)&=\{g\Bb_I\in\Lb_I/\Bb_I\;|\;g^{-1}F(g)\in\Bb_I c_I\Bb_I\}\\
&=\biggl\{(g_1,g_2)\in\GL_i/\Bb_i\times\GL_{d-i}/\Bb_{d-i}\;\big|~g_1^{-1}F(g_1)\in\Bb_ic_1(i)\Bb_i,&\\
&\quad\quad\quad\quad\quad\quad\quad g_2^{-1}F(g_2)\in\Bb_{d-i}c_2(d-i)\Bb_{d-i}\biggl\}\         \ &\\
&\simto  \Omeq^{i-1}\times\Omeq^{d-1-i}
\end{align*}
De m\^eme, on a une description analogue $$Y_{\Lb_I}(\cdo_I)\simto \DL^{i-1}\times\DL^{d-1-i}.$$

\begin{rem}
Soient $I\subset\D$ et $I=\D\ba \{\a_i\}$ pour une racine simple $\a_i\in\D.$ La sous-vari\'et\'e ouverte $X(c_i)$ du diviseur $\o{X(c_i)}$ est une union disjointe de composantes irr\'eductibles isomorphes \`a $X_{\Lb_I}(c_I),$ {\em cf.} \ref{ee}. On notera $C_I$ la composante irr\'eductible fix\'ee par $U_I.$ En vertu du lemme \ref{compact}, on obtient que $\o{C}_I$ est le transform\'e strict de l'espace rationnel $\DG_i$ dans $\o{\O}^{d-1}_{\FM_q}.$
\end{rem}

\subsection{Quotients de Lusztig explicites}\label{R1S3}

Nous calculons explicitement, dans le cas de $\GL_d,$ la construction de Lusztig \cite[2.7, 2.10]{lusztig-coxeter} que nous avons rappel\'ee dans \ref{new9}. En particulier, nous d\'emontrons la propri\'et\'e suivante (voir \ref{bb} pour l'\'enonc\'e plus pr\'ecis):

\begin{thm}
Soient $I=\D\ba \{\a_i\}$ pour une racine simple $\a_i\in\D,$ et $C_I$ la composante irr\'eductible de $X(c_i)$ fix\'ee par $U_I.$ Notons $\Omeq^{d-1,I}:=\Omeq^{d-1}\coprod C_I\subset\o\O^{d-1}_{\FM_q}.$ Alors, il existe un isomorphisme
$$
U_I\ba \Omeq^{d-1,I}\simto \Omeq^{i-1}\times\AM^1\times\Omeq^{d-1-i}.
$$
\end{thm}

\ali\label{a} Soit $\Ub$ le sous groupe de matrices unipotentes triangulaires sup\'erieures de $\GL_d(\oFq),$
$$
\Ub=\biggl\{u=(u_d,u_{d-1},\ldots,u_1)\;\biggl|\; u_d=\left(
                                                \begin{array}{c}
                                                  1 \\
                                                  0\\
                                                  0\\
                                                  0 \\
                                                  0 \\
                                                  \vdots \\
                                                  0 \\

                                                \end{array}
                                              \right)
,~u_{i}=\left(
                                                \begin{array}{c}
                                                  u_{1,i} \\
                                                  \vdots \\
                                                  u_{d-i,i} \\
                                                  1 \\
                                                  0 \\
                                                  \vdots\\
                                                  0\\
                                                \end{array}
                                              \right),
                                              1\leq i\leq d-1\biggl\}.
$$
D'apr\`es la condition \ref{cond de lusztig}, pour que un \'el\'ement $u\in\Ub$ appartienne \`a $X'(c),$ il faut et il suffit qu'il existe $v_i\in\oFq, ~ 1\leq i\leq d-1$ non nuls tels que les conditions suivantes soient v\'erifi\'ees:
$$
\begin{cases}
v_i=u^q_{d-i,i}-u_{d-i,i}\neq0\\
F(u_i)-u_i=v_i\cdot u_{i+1}
\end{cases}
$$
Il est \'evident que pour un tel \'el\'ement $u,$ les $v_i$, $u_i$ et donc $u$ lui-m\^eme sont d\'etermin\'es compl\`etement par sa derni\`ere colonne $u_1.$

\begin{lem}
L'inverse de l'isomorphisme $L$ dans \ref{lusztig thm} (b) est donn\'e par
\begin{align*}
\Omeq^{d-1}&\simto X'(c)\\
\big(x_0:=\frac{X_0}{X_{d-1}},\ldots,x_{d-2}:=\frac{X_{d-2}}{X_{d-1}}\big)&\mapsto u=(u_d,u_{d-1},\ldots,u_1)\text{, o\`u } u_1=\left(
                                                                            \begin{array}{c}
                                                                              x_0 \\
                                                                              \vdots \\
                                                                              x_{d-2} \\
                                                                              1 \\
                                                                            \end{array}
                                                                          \right),\text{ et les $u_i$, $i\geq 2$}\\
                                                                          & \text{sont donn\'es par r\'ecurrence }\begin{cases}u_{i+1}=\frac{1}{v_i}(F(u_i)-u_i)\\
                                                                          v_i=u^q_{d-i,i}-u_{d-i,i}.
                                                                          \end{cases}
\end{align*}
\end{lem}

\begin{preuve}
En sachant que $L:X'(c)\to \O^{d-1}_{\FM_q}$ est un isomorphisme, il suffit de montrer que $\forall u=(u_{ij})\in X'(c),$ $L(u)=(u_{11},\ldots,u_{d-1,1}).$ En effet, $w_\D=\left(
               \begin{array}{ccc}
                  &  & 1 \\
                  & \iddots &  \\
                 1 &  &  \\
               \end{array}
             \right)
,$ donc la premi\`ere colonne de $u\cdot w_\D$ est \'egale \`a $(u_{11},\ldots,u_{d-1,1},1)^t.$ On en d\'eduit l'\'enonc\'e du lemme.
\end{preuve}

Nous pouvons reformuler le th\'eor\`eme \ref{lusztig thm} de la fa\c con suivante:
\begin{coro}
\begin{description}
\item [(a)] Le morphisme
\begin{align*}
X'(c)  \To{} & \GM_m\times\cdots\times\GM_m ~(\text{$d-1$ fois})\\
u=(u_{k,l}) \mapsto & (v_{d-1},\ldots,v_1)
\end{align*}
o\`u $v_i=u_{d-i,i}^q-u_{d-i,i}$ induit un isomorphisme $U\ba X'(c)\simto (\GM_m)^{d-1}.$

\item [(b)] Soient $I\subset\Delta$ et $\Delta\backslash I=\{\alpha_i\}.$ Alors, l'application
\begin{align*}
U_I\ba X'(c) \To{} & \Omeq^{i-1}\times\GM_m\times\Omeq^{d-1-i}\\
u=(u_{k,l})\pmod{U_I} \mapsto & \biggl((u_{1,d+1-i},\ldots,u_{i-1,d+1-i}),v_{d-i}, (u_{i+1,1},\ldots,u_{d-1,1})\biggl)
\end{align*}

est un isomorphisme. En plus, son compos\'e avec l'isomorphisme $\Omeq^{d-1}\simto X'(c)$ induit un isomorphisme
\ini\begin{equation}\label{quotient I}
U_I\ba \Omeq^{d-1}\simto\Omeq^{i-1}\times\GM_m\times\Omeq^{d-1-i}.
\end{equation}
\end{description}
\end{coro}

\ali Soit $I\subset\D$ tel que $\D\ba I=\{\a_i\}.$ On notera $\rho_I$ le morphisme donn\'e par le compos\'e:
$$\Omeq^{d-1}\twoheadrightarrow U_I\ba\Omeq^{d-1}\simto\Omeq^{i-1}\times\GM_m\times\Omeq^{d-1-i}.$$
Rappelons que $C_I$ est la composante irr\'eductible de $X(c_i)$ fix\'ee par $U_I.$ Dans la suite, nous \'etudions $\Omeq^{d-1,I}:=\Omeq^{d-1}\coprod C_I\subset\o{X}(c)$ une compactification partielle de $\Omeq^{d-1}$ associ\'ee \`a $I$, et nous montrons que le morphisme $\rho_I$ se prolonge \`a $\Omeq^{d-1,I}$ induisant un isomorphisme $U_I\ba \Omeq^{d-1,I}\simto\Omeq^{i-1}\times\AM^1\times\Omeq^{d-1-i}$ (voir \ref{bb}). Cette compactification sera utilis\'ee dans la d\'emonstration du th\'eor\`eme \ref{main theorem 2}.


Sous les coordonn\'ees projectives $[X_0: X_1:\ldots:X_{d-1}]$ de $\PM^{d-1}_{\FM_q},$ pour $1\leq j\leq d-1,$ $\Pb_{\D\ba\{\a_j\}}$ est le stabilisateur du sous-espace $\DG_j$ dans $\PM^{d-1}_{\FM_q}.$ Notons que $\o{C}_{\D\ba\{\a_j\}}$ est le transform\'e strict de $\DG_j$ dans $\o{\O}^{d-1}_{\FM_q}.$ Soit $Z$ la r\'eunion des transform\'es stricts de tous les sous-espaces $\FM_q$-rationnels sauf celui de $\DG_j,~\forall 1\leq j\leq d-1.$ Consid\'erons la sous-vari\'et\'e ouverte $\o{\O}^{d-1}_{\FM_q}\ba Z$ de $\o{\O}^{d-1}_{\FM_q}$ qui est en fait la r\'eunion de $\O^{d-1}_{\FM_q}$ avec les composantes $C_{\D\ba\{\a_i\}}$ pour tout $\a_i\in\D.$

\begin{lemme}\label{cc}
La vari\'et\'e $\o{\O}^{d-1}_{\FM_q}\ba Z$ est une sous-vari\'et\'e affine ouverte d'un espace affine $\AM^{d-1}$ de coordonn\'ees affines $(y_1,\ldots,y_{d-1}):=(\frac{X_1}{X_0},\frac{X_2}{X_1},\ldots,\frac{X_{d-1}}{X_{d-2}}).$ De plus, la composante $C_{\D\ba\{\a_i\}}$ est contenue dans l'hyperplan $y_i=0.$
\end{lemme}
\begin{preuve}
Consid\'erons la fonction
$$
\d_d(y_1,\ldots,y_{d-1})=\det\left(
                               \begin{array}{ccccc}
                                 1 & y_1 & y_2 & \cdots  & y_1\cdots y_{d-1} \\
                                 1 & y_1^q & y_2^q & \cdots  & (y_1\cdots y_{d-1})^q \\
                                 \vdots & \vdots & \vdots  & \vdots & \vdots \\
                                 1 & y_1^{q^{d-1}} & y_2^{q^{d-1}} & \cdots  &(y_1\cdots y_{d-1})^{q^{d-1}} \\
                               \end{array}
                             \right)^{q-1}
$$
sur $\AM^{d-1}:=\Spec(\FM_q[y_1,\ldots,y_{d-1}]).$ Le morphisme
\begin{align*}
\varphi:\AM^{d-1}&\To{}\PM^{d-1}_{\FM_q}\\
(y_1,\ldots,y_{d-1})&\longmapsto [1:y_1:y_1y_2:\cdots:y_1\cdots y_{d-1}]
\end{align*}
induit un isomorphisme $\AM^{d-1}_\vide :=\{(y_1,\ldots,y_{d-1})\in\AM^{d-1}~|~ \d_d(y_1,\ldots,y_{d-1})\neq 0\}\simto \O^{d-1}_{\FM_q}.$ Pour chaque $i\in\{1,\ldots,d-1\}$ fix\'e, on a une factorisation
\begin{align*}
\d_d(y_1,\ldots,y_{d-1})&=C\prod_{(a_0,\ldots,a_{d-1})\in\FM^d_q\ba\{0\}}(a_0+a_1y_1+\cdots+a_{d-1}y_1\cdots y_{d-1})\\
&=C\cdot(y_1\cdots y_i)^{q^{d-i}-1}\cdot\d_{d-i}(y_{i+1},\ldots,y_{d-1})\cdot \d_{d,i}(y_1,\ldots,y_{d-1}),
\end{align*}
o\`u $C$ est une constante, et $\d_{d,i}(y_1,\ldots,y_{d-1})$ est le produit de $a_0+a_1y_1+\cdots+a_{d-1}y_1\cdots y_{d-1}$ sur tous les $(a_0,\ldots,a_{d-1})\in\FM^d_q\ba\{0\}$ tels que $a_j\neq 0$ pour certain $j<i.$ Il est clair que $\d_{d,i}\equiv \d_i(y_1,\ldots,y_{i-1})^{q^{d-i}}\pmod{y_i}.$

Consid\'erons $\AM^{d-1,i}$ la sous-vari\'et\'e ouverte de $\AM^{d-1}$ d\'efinie par $$\{(y_1,\ldots,y_{d-1})\in\AM^{d-1}~|~\frac{\d_d(y_1,\ldots,y_{d-1})}{y_i^{q^{d-i}-1}}\neq 0\}.$$
Il r\'esulte de la factorisation de $\d_d$ ci-dessus que $\AM^{d-1,i}=\AM^{d-1}_\vide\coprod A_i,$ o\`u $A_i=\{(y_1,\ldots,y_{d-1})~|~y_i=0,~\d_{d-i}(y_{i+1},\ldots,y_{d-1})\cdot \d_i(y_1,\ldots,y_{i-1})\neq0\}.$ D\'efinissons une fonction $f$ comme le produit de $a_0+a_1y_1+\cdots+a_{d-1}y_1\cdots y_{d-1}$ sur tous les $(a_0,\ldots,a_{d-1})\in\FM^d_q\ba \{0\}$ tels qu'il existe au moins deux $a_j$ non nuls. Alors la sous-vari\'et\'e ouverte $D(f)$ de $\AM^{d-1}$ qui est le compl\'ementaire de $\{f=0\}$ est affine; il s'ensuit que $D(f)=\AM^{d-1,1}\cup\cdots\cup\AM^{d-1,d-1}.$

Consid\'erons ensuite la restriction du morphisme $\varphi$ \`a $D(f),$ not\'ee encore $\varphi:D(f)\to \PM^{d-1}_{\FM_q}.$ Par d\'efinition, l'image r\'eciproque de $\DG_1$ sous $\varphi$ s'identifie \`a $A_1,$ et l'image r\'eciproque des autres points rationnels est vide. Comme $A_1$ est un diviseur de $D(f),$ il r\'esulte de la propri\'et\'e universelle d'\'eclatement qu'il existe un morphisme $\varphi_1:D(f)\to Y_1,$ o\`u $Y_1$ est l'\'eclatement de $\PM^{d-1}_{\FM_q}$ le long de tous ses points rationnels, {\em cf.} \ref{compact}. De plus, $\varphi_1$ induit un isomorphisme entre $A_1$ et le diviseur exceptionnel le long de $\DG_1$ priv\'e les transform\'es stricts de toutes les droites rationnelles sauf celui de $\DG_2.$ De la m\^eme mani\`ere, il existe un morphisme $\varphi_{2}:D(f)\to Y_{2}$ qui induit un isomorphisme entre $A_2$ et le diviseur exceptionnel le long de $\DG_2$ priv\'e le transform\'e strict de tous les plans rationnels sauf celui de $\DG_3.$ Finalement, par r\'ecurrence, on obtient un morphisme $\varphi_{d-2}:D(f)\to Y_{d-2}=\o{\O}^{d-1}_{\FM_q}$ de l'image $\o{\O}^{d-1}_{\FM_q}\ba Z.$ En vertu du \guillemotleft ~main theorem~\guillemotright~ de Zariski \cite[4.4]{Liu}, $\varphi_{d-2}$ induit un isomorphisme $D(f)\simto\o{\O}^{d-1}_{\FM_q}\ba Z.$

\end{preuve}

\begin{rem}
En fait, les coordonn\'ees affines $(y_1,\ldots,y_{d-1})$ comme ci-dessus viennent d'un calcul direct de l'\'eclatement dans le lemme \ref{compact}.
\end{rem}

\begin{lemme}
Sous les coordonn\'ees $(y_1,\ldots,y_{d-1}),$ le morphisme $\rho_I:\O^{d-1}_{\FM_q} \To{}  \Omeq^{i-1}\times\GM_m\times\Omeq^{d-1-i}$ est donn\'e par
\begin{equation*}
(y_1,\ldots,y_{d-1})\longmapsto \biggl(\big(\frac{a_{1,d-i}}{(y_1\cdots y_{i-1})^{q^{d-i}}},\ldots,\frac{a_{i-1,d-i}}{ y_{i-1}^{q^{d-i}}}\big),v_{d-i},\big(\frac{1}{y_{i+1}\cdots y_{d-1}},\ldots,\frac{1}{y_{d-1}}\big)  \biggl),
\end{equation*}
o\`u $a_{jk}\in\FM_q(y_1,\ldots,y_{d-1})$ est donn\'e par r\'ecurrence:
$$
\left\{
             \begin{array}{ll}
               a_{j1}=\frac{1-(y_j\cdots y_{d-1})^{q-1}}{1-y_{d-1}^{q-1}}, & \hbox{$1\leq j\leq d-2$;} \\
               a_{jk}=\frac{a^q_{j,k-1}-(y_j\cdots y_{d-k})^{q^k-q^{k-1}}a_{j,k-1}}{a^q_{d-k,k-1}-y_{d-k}^{q^k-q^{k-1}}a_{d-k,k-1}}, & \hbox{$1\leq j\leq d-k-1$,}
             \end{array}
\right.
$$
et $v_{d-i}=\frac{a^q_{i,d-1-i}-y_{i}^{q^{d-i}-q^{d-1-i}}a_{i,d-1-i}}{y_i^{q^{d-i}}}. $
\end{lemme}

\begin{preuve}
Calculons tout d'abord la matrice unipotente $u\in X'(c)$ dans \ref{a}. \'Ecrivons un point $(x_0,\ldots,x_{d-2})\in \O^{d-1}_{\FM_q}$ sous coordonn\'ees $(y_1,\ldots,y_{d-1})$ de $\AM^{d-1}_{\FM_q}$ ({\em cf.} \ref{cc}),
$$
x_j=\frac{1}{y_{j+1}\cdots y_{d-1}},~0\leq j\leq d-2.
$$
D'apr\`es le lemme \ref{a}, on a
$$u_{j1}=x_{j-1}=\frac{1}{y_j\cdots y_{d-1}},~1\leq j\leq d-1\text{ et }v_1=u^q_{d-1,1}-u_{d-1,1}=\frac{1-y_{d-1}^{q-1}}{y_{d-1}^q}.
$$
On en d\'eduit que, pour $1\leq j\leq d-2,$
$$
u_{j2}=\frac{u^q_{j1}-u_{j1}}{v_1}=\frac{a_{j1}}{(y_j\cdots y_{d-2})^q}, \text{ o\`u $a_{j1}=\frac{1-(y_j\cdots y_{d-1})^{q-1}}{1-y_{d-1}^{q-1}}$}
$$
et
$$
v_2=u^q_{d-2,2}-u_{d-2,2}=\frac{a^q_{d-2,1}-y^{q^2-q}_{d-2}a_{d-2,1}}{y^{q^2}_{d-2}}.
$$
Alors, pour $1\leq j\leq d-3,$ on a
$$
u_{j3}=\frac{u^q_{j2}-u_{j2}}{v_2}=\frac{a_{j2}}{(y_j\cdots y_{d-3})^{q^2}}, \text{ o\`u } a_{j2}=\frac{a^q_{j,1}-(y_j\cdots y_{d-2})^{q^2-q}a_{j1}}{a^q_{d-2,1}-y^{q^2-q}_{d-2}a_{d-2,1}}.
$$
Par r\'ecurrence, posons
$$
a_{jk}=\frac{a^q_{j,k-1}-(y_j\cdots y_{d-k})^{q^k-q^{k-1}}a_{j,k-1}}{a^q_{d-k,k-1}-y_{d-k}^{q^k-q^{k-1}}a_{d-k,k-1}},~1\leq j\leq d-k-1.
$$
On a donc
$$
v_k=\frac{a^q_{d-k,k-1}-y_{d-k}^{q^k-q^{k-1}}a_{d-k,k-1}}{y_{d-k}^{q^k}}\text{ et } u_{j,k+1}=\frac{a_{jk}}{(y_j\cdots y_{d-1-k})^{q^k}}.
$$
\end{preuve}

\begin{lemme}\label{aa}
Pour $j\in\{1,\ldots,i-1\},$ $a_{j,d-i}=\frac{1+f(y_j,\ldots,y_{d-1})}{1+g(y_j,\ldots,y_{d-1})}\in\FM_q(y_j,\ldots,y_{d-1}),$ o\`u $f(y_j,\ldots,y_{d-1})$ et $g(y_j,\ldots,y_{d-1})$ appartiennent \`a $(y_j,\ldots,y_{d-1})\FM_q[y_j,\ldots,y_{d-1}].$ En particulier, $\ord_{y_n}(a_{j,d-i})=0,$ o\`u $\ord_{y_n}$ est la valuation de $\FM_q(y_j,\ldots,y_{d-1})$ d\'efinie par $y_n$ telle que $\ord_{y_n}(y_n)=1,~j\leq n\leq d-1.$ De plus, $a_{j,d-i}\equiv 1\pmod{y_i},$ $v_{d-i}\in\FM_q(y_i,\ldots,y_{d-1}),$ et $\ord_{y_i}(v_{d-i})=-q^{d-i}<0.$
\end{lemme}
\begin{preuve}
D\'emontrons-le par r\'ecurrence. Pour $i=d-1,$ $a_{j1}=\frac{1-(y_j\cdots y_{d-1})^{q-1}}{1-y_{d-1}^{q-1}}$ et $v_1=\frac{1-y_{d-1}^{q-1}}{y_{d-1}^q}.$ Soit $i=k,$ et supposons que l'\'enonc\'e soit vrai pour $i\in\{k+1,\ldots, d-1\}.$ Par d\'efinition, on a $$a_{j,d-k}=\frac{a^q_{j,d-1-k}-(y_j\cdots y_{k})^{q^{d-k}-q^{d-1-k}}a_{j,d-1-k}}{a^q_{k,d-1-k}-y_{k}^{q^{d-k}-q^{d-1-k}}a_{k,d-1-k}}.$$
Notons que par r\'ecurrence, $a_{j,d-1-k}$ et $a_{k,d-1-k}$ appartiennent \`a $\FM_q(y_j,\ldots,y_{d-1}),$ et ils sont de la forme $\frac{1+f}{1+g}$ o\`u $f,g\in (y_j,\ldots,y_{d-1})\FM_q[y_j,\ldots,y_{d-1}].$ On voit imm\'ediatement qu'il existe $\wt{f},\wt{g}\in(y_j,\ldots,y_{d-1})\FM_q[y_j,\ldots,y_{d-1}]$ tels que $a_{j,d-k}=\frac{1+\wt{f}}{1+\wt{g}}\in\FM_q(y_j,\ldots,y_{d-1}),$ donc $\ord_{y_n}(a_{j,d-i})=0,~j\leq n\leq d-1.$ Calculons la valeur $a_{j,d-k}$ modulo $y_k.$ Nous avons
\begin{align*}
a_{j,d-k}&=\frac{a^q_{j,d-1-k}-(y_j\cdots y_{k})^{q^{d-k}-q^{d-1-k}}a_{j,d-1-k}}{a^q_{k,d-1-k}-y_{k}^{q^{d-k}-q^{d-1-k}}a_{k,d-1-k}}\\
&\equiv \big(\frac{a_{j,d-1-k}}{a_{k,d-1-k}}\big)^q\pmod{ y_k}.
\end{align*}
Notons que
\begin{align*}
\frac{a_{j,d-1-k}}{a_{k,d-1-k}}&=\frac{a^q_{j,d-2-k}-(y_j\cdots y_{k+1})^{q^{d-1-k}-q^{d-2-k}}a_{j,d-2-k}}{a^q_{k+1,d-2-k}-y_{k+1}^{q^{d-1-k}-q^{d-2-k}}a_{k+1,d-2-k}}\cdot\frac{a^q_{k+1,d-2-k}-y_{k+1}^{q^{d-1-k}-q^{d-2-k}}a_{k+1,d-2-k}}{a^q_{k,d-2-k}-(y_k y_{k+1})^{q^{d-1-k}-q^{d-2-k}}a_{k,d-2-k}}\\
&\equiv \big(\frac{a_{j,d-2-k}}{a_{k,d-2-k}}\big)^q\pmod{ y_k}\equiv\cdots\equiv\big(\frac{a_{j1}}{a_{i1}}\big)^{q^{d-2-k}}\pmod{y_k}.
\end{align*}
On en d\'eduit que
\begin{align*}
a_{j,d-k}&\equiv \big(\frac{a_{j1}}{a_{i1}}\big)^{q^{d-1-k}}\pmod{y_k}\\
&= \biggl(\frac{1-(y_j\cdots y_{d-1})^{q-1}}{1-y_{d-1}^{q-1}}\cdot\frac{1-y_{d-1}^{q-1}}{1-(y_k\cdots y_{d-1})^{q-1}}\biggl)^{q^{d-1-k}}\\
&\equiv 1\pmod{y_k}.
\end{align*}

Rappelons que $$v_{d-k}=\frac{a^q_{k,d-1-k}-y_{k}^{q^{d-k}-q^{d-1-k}}a_{k,d-1-k}}{y_k^{q^{d-k}}}\in\FM_q(y_k,\ldots,y_{d-1}).$$
Par r\'ecurrence, on sait que $\ord_{y_k}(a_{k,d-1-k})=0.$ On en d\'eduit que $\ord_{y_k}(v_{d-k})=-q^{d-k}.$ Ceci termine la preuve du lemme.

\end{preuve}

\begin{coro}\label{bb}
Soit $I=\D\ba\{\a_i\},$ consid\'erons le morphisme
\begin{align*}
\rho'_I:\O^{d-1}_{\FM_q}& \To{}  \Omeq^{i-1}\times\GM_m\times\Omeq^{d-1-i}\\
(y_1,\ldots,y_{d-1})&\longmapsto \biggl(\big(\frac{a_{1,d-i}}{(y_1\cdots y_{i-1})^{q^{d-i}}},\ldots,\frac{a_{i-1,d-i}}{ y_{i-1}^{q^{d-i}}}\big),v_{d-i}^{-1},\big(\frac{1}{y_{i+1}\cdots y_{d-1}},\ldots,\frac{1}{y_{d-1}}\big)  \biggl)
\end{align*}
qui est le compos\'e de $\rho_I$ avec le morphisme $\id\times()^{-1}\times\id:\Omeq^{i-1}\times\GM_m\times\Omeq^{d-1-i}\to \Omeq^{i-1}\times\GM_m\times\Omeq^{d-1-i}.$ Avec la m\^eme formule, $\rho'_I$ se prolonge en un morphisme $$\o{\rho}_I:\O^{d-1,I}_{\FM_q}\to \Omeq^{i-1}\times\AM^1\times\Omeq^{d-1-i}.$$ De plus,
\begin{description}
  \item[(a)] La restriction de $\o{\rho}_I$ \`a $C_I$ est un morphisme radiciel ({\em cf.} \cite[AG. 18.2]{borel}), donn\'e par
\begin{align*}
C_I=\O^{i-1}_{\FM_q}\times\{0\}\times\O^{d-1-i}_{\FM_q} \To{} & \O^{i-1}_{\FM_q}\times\{0\}\times\O^{d-1-i}_{\FM_q}\\
\Big((y_1,\ldots, y_{i-1}),0,(y_{i+1},\ldots,y_{d-1})\Big)   \longmapsto & \Big(\big((\frac{1}{y_1\cdots y_{i-1}})^{q^{d-i}},(\frac{1}{y_2\cdots y_{i-1}})^{q^{d-i}},\ldots,(\frac{1}{y_{i-1}})^{q^{d-i}}\big), \\
&  0, \big(\frac{1}{y_{i+1}\cdots y_{d-1}},\frac{1}{y_{i+2}\cdots y_{d-1}},\ldots,\frac{1}{y_{d-1}} \big) \Big).
\end{align*}
  \item[(b)] $\o{\rho}_I$ est surjectif, et induit un isomorphisme $U_I\ba \O^{d-1,I}_{\FM_q}\simto \Omeq^{i-1}\times\AM^1\times\Omeq^{d-1-i},$ o\`u $U_I\ba \O^{d-1,I}_{\FM_q}$ est le quotient de la $\FM_q$-vari\'et\'e affine $\O^{d-1,I}_{\FM_q}$ par le groupe fini $U_I.$
\end{description}
$$\xymatrix{
\O^{d-1}_{\FM_q} \ar@{^(->}[r] \ar[d]^{\rho_{I}}_{/U_I} & \Omeq^{d-1,I} \ar@{<-^)}[r] \ar@{.>}[d]^{\o{\rho}_I} & C_I=\O^{i-1}_{\FM_q}\times\{0\}\times\O^{d-1-i}_{\FM_q} \ar@{.>}[d]\\
\O^{i\!-\!1}_{\FM_q}\!\times\!\GM_m\!\times\!\O^{d\!-\!1\!-\!i}_{\FM_q} \ar@{^(->}[r]^{\id\!\times\!()^{\!-\!1}\!\times\!\id}  & \O^{i\!-\!1}_{\FM_q}\!\times\!\AM^1\!\times\!\O^{d\!-\!1\!-\!i}_{\FM_q} \ar@{<-^)}[r]^{\id} & \O^{i-1}_{\FM_q}\times\{0\}\times\O^{d-1-i}_{\FM_q}
}$$
\end{coro}
\begin{preuve}
(a) D'apr\`es le lemme \ref{aa}, $a_{j,d-i}\equiv 1\pmod{y_i},~1\leq j\leq i-1$ et $v^{-1}_{d-i}\equiv 0\pmod{y_i}.$ Comme la composante $C_I$ est contenue dans l'hyperplan $y_i=0$ (voir \ref{cc}), $\rho'_I$ se prolonge \`a $\O^{d-1,I}_{\FM_q}$ de la mani\`ere \'enonc\'ee comme ci-dessus.

(b) Comme $\O^{d-1,I}_{\FM_q}$ est une vari\'et\'e affine normale et $U_I$ est fini, le quotient $U_I\ba \O^{d-1,I}_{\FM_q}$ existe ({\em cf.} \cite[I 6.15]{borel}) et il est normal. Notons que le morphisme $\o{\rho}_I$ est constant sur les $U_I$-orbites, il se factorise donc par le quotient $U_I\ba \O^{d-1,I}_{\FM_q},$ et induit un morphisme radiciel birationnel $\a:U_I\ba \O^{d-1,I}_{\FM_q}\onto \Omeq^{i-1}\times\AM^1\times\Omeq^{d-1-i};$ il s'ensuit que $\a$ est un isomorphisme en vertu du \guillemotleft ~main theorem~\guillemotright~ de Zariski \cite[4.4]{Liu}.

\end{preuve}

\section{Th\'eor\`eme principal}

\subsection{\'Enonc\'e du th\'eor\`eme}

Pour all\'eger un peu les notations, on notera $\O :=\O^{d-1}_{\FM_q}$ et $\DL :=\DL^{d-1}.$ Consid\'erons la compactification $\o{\O}$ de $\O$ et sa normalisation $\o{\DL}$ dans $\DL$ ({\em cf.} \cite{bonnafe-rouquier-compact} ou \ref{new5}). Soient $I$ un sous-ensemble propre de $\D$ et $C_I$ la composante irr\'eductible de $X(\prod_{i\in I}s_{\a_i})\subset\o{\O}^{d-1}_{\oFq}$ fix\'ee par $U_I,$ o\`u le produit est pris suivant l'ordre $s_{\a_1}, s_{\a_2},\ldots,s_{\a_{d-1}}.$ Rappelons que l'on a un diagramme commutatif:
$$\xymatrix{
\DL\ar[d]_\pi \ar@{^(->}[r]^{j'} &\o{\DL} \ar@{<-^)}[r]^{j'_I} \ar[d]_{\o{\pi}} & \o{\pi}^{-1}(C_I) \ar[d]_{\o{\pi}_I}\\
\O \ar@{^(->}[r]^j &\o{\O} \ar@{<-^)}[r]^{j_I} & C_I
}$$
o\`u $\o{\pi}^{-1}(C_I):=(\o{\DL}\times_{\o{\O}} C_I)_{\re}.$ Dans cette partie, on se consacre \`a d\'emontrer la conjecture A' dans ce cas.

\begin{theo}\label{main theorem 2}
Le morphisme de restriction
$$
R\G(\DL,\L)\To{\res.}R\G(\o\pi^{-1}(C_I),j'^*_IRj'_*\L)
$$
induit un isomorphisme
$$
R\G(\DL,\L)^{U_I}\simto R\G(\o\pi^{-1}(C_I),j'^*_IRj'_*\L).
$$
\end{theo}

Comme un sous-groupe de Levi rationnel de $\GL_d$ est un produit de $\GL_n,$ d'apr\`es \ref{R1L2} et \ref{R1S5}, il suffit d'\'etablir le cas de codimension $1.$

\subsection{\'Etape 1: le cas de codimension 1}\label{ff}

On reprend l'id\'ee de \cite{bonnafe-rouquier-coxeter} et de Dudas \cite{dudas}. Lorsque $C_I$ est une composante de codimension $1$, il existe une racine $\a_i$ telle que $I=\D\ba\{\a_i\}.$ Notons $\O^I$ (resp. $\DL^I$) la compactification partielle $\O\coprod C_I$ (resp. $\DL\coprod \o{\pi}^{-1}(C_I)$). Consid\'erons la pr\'eimage $Y^0$ d'une composante connexe de $U\ba \DL$ ({\em cf.} \cite[3.2]{bonnafe-rouquier-coxeter} et \cite[Prop. 4.53]{dudas}). Ici $Y^0$ est une composante connexe de $\DL,$ finie \'etale au-dessus de $\O$ de groupe de Galois $$H:=\Ker(\norm:\FM_{q^d}^\times\to\FM_q^\times)\cong\ZM/(1+q+\cdots+q^{d-1}).$$ Le tore non-d\'eploy\'e $T_d\simto\FM_{q^d}^\times$ ({\em cf.} section \ref{new8}) agit transitivement sur les composantes connexes de $U\ba \DL.$ En notant que $H$ est le stabilisateur de la composante $U\ba Y^0$ dans $T_d$ (c'est aussi le stabilisateur de $Y^0$), on obtient un isomorphisme $\GL_d(\FM_q)\times(T_d\rtimes \langle F\rangle)^{\opp}$-\'equivariant (ici $F$ d\'esigne l'endomorphisme de Frobenius):
$$Y^0\times_H T_d\simto \DL. $$
Rappelons que le groupe fondamental mod\'er\'e de $\GM_m^{d-1}$ est le produit $(d-1)$ fois du groupe fondamental mod\'er\'e de $\GM_m.$ Alors, il existe $d-1$ entiers positifs $m_1,\ldots,m_{d-1}$ divisant $|H|,$ et un rev\^etement $\varpi:\GM_m^{d-1}\onto\GM_m^{d-1}$ de groupe de Galois $\prod_{j}\mu_{m_j}$ tels que le rev\^etement galoisien $U\ba Y^0\onto U\ba \O=\GM_m^{d-1}$ soit un quotient de $\varpi.$ On note $N$ le groupe de Galois du rev\^etement $\GM_m^{d-1}\onto U\ba Y^0$ et on a un isomorphisme canonique $(\prod_j\mu_{m_j})/N\simto H.$ En diminuant les $m_j,$ on peut supposer et on le fera que la restriction $\phi_j:\mu_{m_j}\to H$ du morphisme canonique $\prod_j\mu_{m_j}\onto H$ est injective pour tout $j.$ On obtient un diagramme commutatif:

$$\xymatrix{\GM_m^{d-1}\ar@{>>}[dr]_{/\prod_{j}\mu_{m_j}}^{\varpi}\ar@{->>}[r]^{/N} & U\ba Y^0\ar@{>>}[d]^{/H}  \\
& \GM_m^{d-1} =U\ba \O }
$$

Notons $Y^1$ le produit fibr\'e de $\GM_m^{d-1}$ et $\O$ au-dessus de $\GM_m^{d-1}$ via $\varpi$ et l'application de quotient par $U.$ On peut alors former le diagramme suivant, dans lequel les carr\'es sont cart\'esiens:
$$\xymatrix{
Y^1\ar@{>>}[r]\ar@{>>}[d]^{\pi^1}\ar@(l,l)[dd]^{\pi_1}  & U_I\ba Y^1\ar@{>>}[r]\ar@{>>}[d] & \GM^{d-1}_m \ar@{>>}[d]\ar@(ru,u)@{>>}^{/\prod\mu_{m_j}}[rdd] \\
Y^0\ar@{>>}[r]\ar@{>>}[d]^{\pi_0}  &U_I\ba Y^0\ar@{>>}[r]\ar@{>>}[d] & U\ba Y^0 \ar@{>>}[d] \\
\O\ar@{>>}[r]  &U_I\ba\O\ar@{>>}[r] & U\ba \O \ar[r]^\thicksim &\GM_m^{d-1}
}
$$

Via l'identification $U_I\ba \O\simto (\O^{i-1}\times\O^{d-1-i})\times\GM_m$ ({\em cf.} \ref{quotient I}), l'application $U_I\ba \O\onto V_I\ba (U_I\ba\O)=\GM_m^{d-1}$ se d\'ecompose en $\pi_{V_I}\times\id_{\GM_m}$ o\`u $V_I$ est le radical unipotent de $L_I$ (voir section \ref{new9}), $\pi_{V_I}:\O^{i-1}\times\O^{d-1-i}\onto \GM_m^{d-2}$ est le morphisme de quotient par $V_I.$ L'avantage de cette construction est que si l'on forme le produit fibr\'e suivant:
$$\xymatrix{Y^1_I\ar@{>>}[r]\ar@{>>}[d] & \GM_m^{d-2}\ar@{>>}[d]^{/\prod_{j\neq i}\mu_{m_j}}\\
\O^{i-1}\times\O^{d-1-i}\ar@{>>}[r] & \GM_m^{d-2}
}
$$
on peut d\'ecomposer le quotient $U_I\ba Y^1\simto Y^1_I\times\GM_m$ de mani\`ere compatible avec la d\'ecomposition $U_I\ba \O=(\O^{i-1}\times\O^{d-1-i})\times \GM_m.$ Par construction, $Y^1_I$ est \'etale au-dessus de la vari\'et\'e normale $\O^{i-1}\times\O^{d-1-i},$ et donc elle est normale.

On peut alors former un diagramme commutatif:
$$\xymatrix
{U_I\ba\! Y^1\!=\!Y^1_I\!\times\!\GM_m\!\ar@{^(->}[r]^{   j'_1}\ar[d] & Y^1_I\times\AM^1\ar@{<-^)}[r]^{i'_{1,I}}\ar[dd] & Y^1_I\times\{0\}\ar[dd]\\
U_I\ba Y^0\ar[d]\\
\O^{i-1}\! \times\!\GM_m\!\times\!\O^{d-1-i} \ar@{^(->}[r] & \O^{i-1}\times\AM^1\times\O^{d-1-i} \ar@{<-^)}[r] & \O^{i-1}\times\{0\}\times\O^{d-1-i}
}$$
dont le morphisme $\AM^1\to\AM^1$ au milieu est donn\'e par $x\mapsto x^{m_i}.$ La vari\'et\'e normale $Y^1_I\times\AM^1$ s'identifie alors \`a la normalisation de $\O^{i-1}\times\AM^1\times\O^{d-1-i}$ dans $U_I\ba Y^1.$

Notons $\o{Y^0}\To{\o{\pi}_0} \O^I$ la normalisation de $\O^I$ dans $Y^0$ et $\o{\pi}^{-1}_0(C_I)$ le sous sch\'ema ferm\'e avec la structure r\'eduite de $\o{Y^0}.$ Comme $\O^I$ est affine (voir le lemme \ref{cc}), $\o{Y^0}$ est affine. Comme $Y^0$ est une composante connexe de $\DL,$ et $(\DL^I,\o{\pi}^{-1}(C_I))$ est un couple lisse ({\em cf.} \cite[Thm. 1.2 (a)]{bonnafe-rouquier-compact}), on sait alors que $(\o{Y^0},\o{\pi}^{-1}_0(C_I))$ est un couple lisse. Ceci nous fournit un diagramme commutatif:
$$\xymatrix{
Y^0\ar@{^(->}[r]^{j_0}\ar[d]^{\pi_0} & \o{Y^0}\ar@{<-^)}[r]^{j_{0,I}}\ar[d]^{\o{\pi}_0} & \o{\pi}^{-1}_0(C_I)\ar[d]\\
\O \ar@{^(->}[r]^{j} & \O^I \ar@{<-^)}[r]^{j_{I}} & C_I
}$$
En notant que $U_I$ est un groupe fini, on peut former le quotient $U_I\ba \o{Y^0}$ qui est une vari\'et\'e normale, et elle s'identifie donc \`a la normalisation de $U_I\ba \O^I=\O^{i-1}\times\AM^1\times\O^{d-1-i}$ ({\em cf.} Corollaire \ref{bb} (b)) dans $U_I\ba Y^0.$ En plus, $U_I$ agit trivialement sur $\o{\pi}^{-1}_0(C_I),$ 
 donc le quotient $\o{\pi}^{-1}_0(C_I)\onto U_I\ba\o{\pi}^{-1}_0(C_I)$ est un morphisme radiciel. On a alors un diagramme commutatif:
$$\tiny\xymatrix
{U_I\ba\! Y^1\!=\!Y^1_I\!\times\!\GM_m\!\ar@{^(->}[r]^{ j'_1}\ar@(ld,lu)[dd]_{/\prod_{j}\mu_{m_j}}\ar[d]^{/N} & Y^1_I\times\AM^1\ar@{<-^)}[r]^{i'_{1,I}}\ar[d] & Y^1_I\times\{0\}\ar[d]_{/N'}\ar@(rd,ru)[dd]^{~/\prod_{j\neq i}\mu_{m_j}}\\
U_I\ba Y^0\ar[d]^{/H}\ar@{^(->}[r] & U_I\ba \o{Y^0}\ar@{<-^)}[r]\ar[d] & U_I\!\ba\!\o{\pi}^{-1}_0\!(\!C_I\!)\ar[d]_{/H'}\\
\O^{i-1}\! \times\!\GM_m\!\times\!\O^{d-1-i} \ar@{^(->}[r] & \O^{i-1}\times\AM^1\times\O^{d-1-i} \ar@{<-^)}[r] & \O^{\!i\!-\!1\!}\!\times\!\{\!0\!\}\!\times\!\O^{\!d\!-\!1\!-\!i\!}
}$$

\begin{lem}
Le morphisme $Y^1_I\times\AM^1\to U_I\ba \o{Y^0}$ est \'etale.
\end{lem}
\begin{preuve}
Notons $N'$ le groupe de Galois de $Y^1_I\times\{0\}$ au-dessus de $U_I\ba\o{\pi}^{-1}_0(C_I)$ et $H'$ le quotient de $\prod_{j\neq i}\mu_{m_j}$ par $N'.$ Par construction, $H$ agit sur $\o{Y^0}$ commutant avec l'action de $U_I,$ donc $H$ agit sur le quotient $U_I\ba\o{Y^0}.$ D'apr\`es \cite[Thm. 1.2 (c)]{bonnafe-rouquier-compact}, le stabilisateur dans $H$ d'un \'el\'ement de $U_I\ba\o{\pi}^{-1}_0(C_I)$ est \'egal \`a $N_c(Y_{c,c_I})\cap H,$ o\`u $N_c(Y_{c,c_I})$ est un sous-groupe de $T_d$ d\'efini dans {\em loc. cit.} section 1.3. En fait, $N_c(Y_{c,c_I})\subset H$ ({\em cf.} \cite[Prop. 3.5]{bonnafe-rouquier-coxeter}), on en d\'eduit que $H'=H/N_c(Y_{c,c_I}).$ En d'autres termes, on a un diagramme dont la deuxi\`eme et la troisi\`eme lignes et la deuxi\`eme et la troisi\`eme colonnes sont exactes.
$$\xymatrix{
 & & 1\ar[d] & 1\ar[d] &\\
 & & \mu_{m_i} \ar[d] \ar[r]^{\phi_i} & N_c(Y_{c,c_I}) \ar[d] &\\
1\ar[r] &N\ar[d] \ar[r] & \prod_{j}\mu_{m_j}\ar[d] \ar[r] & H\ar[d]\ar[r] & 1\\
1 \ar[r] & N'\ar[r] & \prod_{j\neq i}\mu_{m_j}\ar[d] \ar[r] & H'=H/N_c(Y_{c,c_I}) \ar[d]\ar[r] &1\\
& & 1 & 1
}$$
Rappelons que $\phi_i$ est injective. D'apr\`es \cite[Prop. 3.5]{bonnafe-rouquier-coxeter}, $\phi_i$ est aussi sujective. Donc $\phi_i$ est un isomorphisme. Ceci implique que $N\simto N'.$ D'apr\`es \cite[Exp. V Prop. 2.6]{SGA1}, $Y^1_I\times\AM^1\to U_I\ba \o{Y^0}$ est un morphisme \'etale.
\end{preuve}


Posons $\o{Y^1}:=\o{Y^0}\times_{U_I\ba\o{Y^0}}(Y^1_I\times\AM^1)$ le produit fibr\'e de $\o{Y^0}$ et $Y^1_I\times\AM^1$ au-dessus de $U_I\ba\o{Y^0},$ et $\o{\pi}^1:\o{Y^1}\to\o{Y^0}$ la projection vers $\o{Y^0}.$ Alors $\o{\pi}^1$ est un morphisme \'etale. Notons $\o{\pi}_1:=\o{\pi}_0\circ\o{\pi}^1,$ et $\o{\pi}^{-1}_1(C_I)$ le sous-sch\'ema ferm\'e avec la structure r\'eduite de $\o{Y^1}.$
\begin{lemme}\label{lisse}
$(\o{Y^1},\o{\pi}^{-1}_1(C_I))$ est un couple lisse.
\end{lemme}
\begin{preuve}
D'apr\`es le lemme pr\'ec\'edent, $\o{Y^1}$ est \'etale au-dessus de $\o{Y^0}.$ On d\'eduit l'\'enonc\'e du lemme du fait que $(\o{Y^0},\o{\pi}^{-1}_0(C_I))$ est un couple lisse.
\end{preuve}

Consid\'erons le diagramme suivant:
$$
\xymatrix{
                 & Y^1 \ar[rr]^{j_1} \ar[ld]^{\pi^1}\ar[dd]  &           &\o{Y^1} \ar[ld]^{\o{\pi}^1} \ar[dd]   \\
Y^0\ar@{^(->}[rr]^{j_0} \ar[dd]     &              & \o{Y^0} \ar[dd]                \\
                 &U_I\ba Y^1\ar@{^(->}[rr]^{j'_1} \ar[ld]    &    & Y^1_I\times\AM^1\ar[ld]\\
U_I\ba Y^0\ar@^{^(->}[rr]        &        & U_I\ba\o{Y^0}
}
$$
o\`u le morphisme $j_1:Y^1\to\o{Y^1}$ est donn\'e par le produit des morphismes $Y^1\to Y^0\to\o{Y^0}$ et $Y^1\to U_I\ba Y^1\to Y^1_I\times\AM^1.$

\begin{lem}
Les carr\'es
$$
\xymatrix{
Y^1 \ar[r]^{j_1}\ar[d] & \o{Y^1}\ar[d]  &\text{et}  & Y^1 \ar[r]^{j_1}\ar[d]    & \o{Y^1}\ar[d]\\
Y^0 \ar[r]       &\o{Y^0}  &   & U_I\ba Y^1 \ar[r]  & Y^1_I\times\AM^1
}
$$
sont cart\'esiens
\end{lem}
\begin{preuve}
Tout d'abord, notons que le carr\'e
$$
\xymatrix{
Y^0\ar[r]\ar[d] & \o{Y^0}\ar[d]\\
U_I\ba Y^0\ar[r] & U_I\ba \o{Y^0}}
$$
est cart\'esien, car $Y^0$ est l'image r\'eciproque de l'ouvert $U_I\ba Y^0$ dans $\o{Y^0}.$ Pour le premier carr\'e, on a donc
\begin{align*}
\o{Y^1}\times_{\o{Y^0}}Y^0&=(Y^1_I\times\AM^1)\times_{U_I\ba \o{Y^0}}\o{Y^0}\times_{\o{Y^0}}Y^0\\
&=(Y^1_I\times\AM^1)\times_{U_I\ba \o{Y^0}}(U_I\ba Y^0)\times_{U_I\ba Y^0}Y^0\\
&=U_I\ba Y^1\times_{U_I\ba Y^0}Y^0=Y^1.
\end{align*}
Pour le deuxi\`eme carr\'e, on a
\begin{align*}
\o{Y^1}\times_{Y^1_I\times\AM^1}U_I\ba Y^1&=\o{Y^1}\times_{U_I\ba \o{Y^0}}(Y^1_I\times\AM^1)\times_{Y^1_I\times\AM^1}U_I\ba Y^1\\
&=\o{Y^0}\times_{U_I\ba \o{Y^0}}U_I\ba Y^1\\
&=\o{Y^0}\times_{U_I\ba \o{Y^0}}U_I\ba Y^0\times_{U_I\ba Y^0}U_I\ba Y^1\\
&=Y^0\times_{U_I\ba Y^0}U_I\ba Y^1=Y^1.
\end{align*}
D'o\`u l'\'enonc\'e du lemme.
\end{preuve}



D'apr\`es ce qui pr\'ec\`ede, on obtient le diagramme suivant:

$$\tiny\xymatrix @C=.5mm@R=3mm{
&\DL\ar[lddd]_{\pi} \ar@{^(->}[rrr]^{j'} & & & \DL^I\ar[lddd]_{\o{\pi}} \ar@{<-^)}[rrr]^{j'_I} & & & \o{\pi}^{-\!1}(\!C_I\!)\ar[lddd]_{\o{\pi}_I}   \\
& & Y^1\!\ar@(lu,ru)[lldd]^{\pi_1}\ar[ddd]_{\rho_1}\ar[ld]^{\pi^1\; /N} \ar@{^(->}[rrr]^{j_1} & & &  \o{Y^1}\!\ar[ld]\ar[ddd]_{\o{\rho}_1}\ar@(lu,ru)[lldd]^{\o{\pi}_1} \ar@{<-^)}[rrr]^{j_{1,I}}  & &  & \o{\pi}^{-1}_{\mathclap{1}}(\!C_I\!)\ar[ddd]_{\o{\rho}_{1,I}}\ar[ld]^{\o{\pi}^1_I} \\
& Y^0\ar@{^(->}[uu]\ar[ddd]\ar[ld]_{\pi_0}^{/H}\ar@{^(->}[rrr] ^{j_0} & & &  \o{Y^0}\ar[ddd]\ar[ld]^{\o{\pi}_0}\ar@{^(->}[uu] \ar@{<-^)}[rrr]^{j_{0,I}} &  & & \o{\pi}^{-1}_{\mathclap{0}}\!(\!C_I\!)\ar[ld] \ar[ddd]\ar@{^(->}[uu] & &  \\
\O\ar[ddd]_{/ U_I}^{\rho_I} \ar@{^(->}[rrr]^j & & & \O^I\ar[ddd]^{\o{\rho}_I}\ar@{<-^)}[rrr]^{j_I} & & & C_I\ar[ddd]^{\o{\rho}_I|_{C_I}}   &\\
& & U_I\ba Y^1\ar[ddd]\ar[ld]\ar@{^(->}[rrr]^{j'_1}& & &  Y^1_I\!\!\times\!\!\AM^1\ar[ddd]\ar[ld]\ar@{<-^)}[rrr]^{i'_{1,I}}  & &  & Y^1_I\!\!\times\!\!\{\!0\!\}\ar[ld]\ar[ddd] \\
& U_I\!\ba\! Y^0 \ar[ld]\ar[ddd]\ar[rrr] & & & U_I\!\ba\! \o{Y^0}\ar[ld] \ar@{<-^)}[rrr]& & & U_I\!\ba\!\o{\pi}^{-1}_{\mathclap{0}}(\!C_I\!)\ar[ld] &  & \\
\O^{i\!-\!1}\!\!\times\!\!\GM_m\!\!\times\!\!\O^{d\!-\!1\!-\!i}\ar[ddd]_{/ V_I}\ar@{^(->}[rrr] &  & &  \O^{i\!-\!1}\!\!\times\!\!\AM^1\!\!\times\!\!\O^{d\!-\!1\!-\!i}\ar[ddd] \ar@{<-^)}[rrr] & & & \O^{i\!-\!1}\!\!\times\!\!\{\!0\!\}\!\!\times\!\!\O^{d\!-\!1\!-\!i}\ar[ddd] \\
& & \GM_m^{d\!-\!1}\ar@{^(->}[rrr]\ar[ld]_{/N}\ar@(lu,u)[lldd]^{/ \prod\mu_{m_j}} & & &  \GM_m^{i\!-\!1}\!\!\times\!\AM^1\!\times\!\!\GM_m^{d\!-\!1\!-\!i}\ar[lldd] \ar@{<-^)}[rrr] & & &  \GM_m^{i\!-\!1}\!\!\times\!\!\{\!0\!\}\!\!\times\!\!\GM_m^{d\!-\!1\!-\!i}\ar[lldd]  \\
& U\ba Y^0\ar[ld] & & & & & & & & \\
\GM_m^{d\!-\!1}\ar@{^(->}[rrr] & & & \GM_m^{i\!-\!1}\!\!\times\!\AM^1\!\times\!\!\GM_m^{d\!-\!1\!-\!i} \ar@{<-^)}[rrr] & & & \GM_m^{i\!-\!1}\!\!\times\!\!\{\!0\!\}\!\!\times\!\!\GM_m^{d\!-\!1\!-\!i}
}
$$

\ali On commence maintenant \`a d\'emontrer le cas de codimension $1.$ Consid\'erons le diagramme commutatif suivant:
$$\xymatrix{Y^1\ar@{^(->}[r]^{j_1}\ar[d]^{\rho_1} & \o{Y^1}\ar@{<-^)}[r]^{j_{1,I}}\ar[d]^{\o{\rho}_1} & \o{\pi}^{-1}_1(C_I)\ar[d]^{\o{\rho}_{1,I}}\\
U_I\ba\! Y^1\!=\!Y^1_I\!\times\!\GM_m\! \ar@{^(->}[r]^{~j'_1} & Y^1_I\times\AM^1 \ar@{<-^)}[r]^{i'_{1,I}} & Y^1_I\times\{0\}
}$$
Comme $\o{\rho}_1$ est fini, le morphisme canonique $\L\to\rho_{1,*}\rho^*_1\L$ des faisceaux \'etales sur $U_I\ba Y^1$ induit un diagramme commutatif:
$$\xymatrix{R\G(Y^1,\L)\ar[r]^{\res .}_{(4)}  & R\G(\o{\pi}^{-1}_1(C_I),j^*_{1,I}Rj_{1 *}\L)\\
R\G(Y^1_I\times\GM_m,\L)\ar[r]^{\res .}_{(2)}\ar[u]_{(1)} & R\G(Y^1_I\times\{0\},i'^*_{1,I}Rj'_{1 *}\L)\ar[u]_{(3)}
}$$

\begin{lem}
Le morphisme $(4)$ induit un isomorphisme:
\begin{equation*}
R\G(Y^1,\L)^{U_I}\simto R\G(\o{\pi}^{-1}_1(C_I),j^*_{1,I}Rj_{1,I *}\L).
\end{equation*}
\end{lem}
\begin{preuve}
\begin{itemize}
\item\'Etape 1: Notons tout d'abord que l'on a un isomorphisme de faisceaux \'etales sur $U_I\ba Y^1$ $$\L\simto (\rho_{1,*}\L)^{U_I},$$ o\`u $\rho_1$ est fini \'etale de groupe de Galois $U_I$ avec $|U_I|$ inversible dans $\L.$ Donc $(1)$ fournit un isomorphisme $$R\G(Y^1_I\times\GM_m,\L)\simto R\G(Y^1,\L)^{U_I}.$$
\item\'Etape 2: Notons que le diagramme:
$$\xymatrix{
Y^1_I\times\GM_m\ar[r]^{j'_1=\id\times j_2} & Y^1_I\times\AM^1 & Y^1_I\times\{0\} \ar[l]_{i'_{1,I}=\id\times j_2}
}$$
est un changement de base des inclusions:
$$\xymatrix{
\GM_m\ar[r]^{j_2} & \AM^1 & \{0\} \ar[l]_{i_2}
}$$
qui induisent un isomorphisme canonique: $R\G(\GM_m,\L)\simto R\G(\{0\},i_2^*Rj_{2*}\L).$ Consid\'erons ensuite le diagramme suivant dont les carr\'es sont cart\'esiens:
$$\xymatrix{
Y^1_I \ar[r]^{\id} & Y^1_I  & Y^1_I \ar[l]_{\id}  \\
Y^1_I\times\GM_m \ar[u]^{p_1}\ar[d]_{p_2} \ar[r]^{\id\times j_2}  & Y^1_I\times\AM^1 \ar[u]^{p'_1}\ar[d]_{p'_2}  & Y^1_I\times\{0\} \ar[u]^{p''_1}\ar[d]_{p''_2} \ar[l]_{\id\times i_2}  \\
\GM_m \ar[r]^{j_2} & \AM^1  & \{0\} \ar[l]_{i_2}.
}$$
Par la formule de K\"unneth, on a
\begin{align*}
{R\G(Y^1_I\times\GM_m,\L)}&= R\G(Y^1_I,\L)\otimes^L R\G(\GM_m,\L)= R\G(Y^1_I,\L)\otimes^L R\G(\{0\},i_2^*R j_{2*}\L)\\
&= R\G(Y^1_I\times \{0\}, p''^*_1\L\otimes p''^*_2 i_2^*R j_{2*}\L)\\
&= R\G(Y^1_I\times \{0\}, (\id\times i_2)^*(p'^*_1\L\otimes p'^*_2 R j_{2*}\L))  \\
(\text{$p'_2$ est lisse.})&= R\G(Y^1_I\times \{0\}, (\id\times i_2)^*(p'^*_1\L\otimes R(\id\times j_2)_* p_2^*\L))\\
(\text{formule de projection}) &= R\G(Y^1_I\times \{0\},(\id\times i_2)^*R(\id\times j_2)_*((\id\times j_2)^*p'^*_1\L\otimes p_2^*\L))\\
&=R\G(Y^1_I\times\{0\},i'^*_{1,I}Rj'_{1 *}\L)
\end{align*}
i.e. $(2)$ est un isomorphisme.

\item\'Etape 3: Notons que $U_I$ est un $p$-groupe fini. Posons $e_1\in\L[U_I]$ l'idempotent central associ\'e \`a la repr\'esentation triviale de $U_I,$ et $e'_1:=1-e_1\in\L[U_I].$ Le but est de montrer que
$$
\o{\rho}_{1,I,*}j_{1,I}^*Rj_{1*}\L=i'^*_{1,I}Rj'_{1*}\L.
$$
Tout d'abord, notons $D(?,\L-U_I)$ la cat\'egorie d\'eriv\'ee des faisceaux \'etales de $\L$-modules $U_I$-\'equivariants sur $?.$ Le foncteur
    $$
   j_{1,I}^*Rj_{1*}: D^+(Y^1,\L-U_I)\To{} D^+(\o{\pi}^{-1}_1(C_I),\L-U_I)
    $$
est $U_I$-\'equivariant, et induit une d\'ecomposition:
$$
j_{1,I}^*Rj_{1*}\L=e_1(j_{1,I}^*Rj_{1*}\L)\oplus e'_1(j_{1,I}^*Rj_{1*}\L).
$$
D'apr\`es le th\'eor\`eme de puret\'e relative (\cite[Exp. XVI]{SGA4-3}) par rapport au couple lisse $(\o{Y^1},\o{\pi}^{-1}_1(C_I))$ ({\em cf.} le lemme \ref{lisse}), on a des isomorphismes de $\L$-modules:
\begin{align*}
\a:&\L\simto j_{1,I}^*R^0j_{1*}\L,\\
\b:&\L(-1)\simto j_{1,I}^*R^1j_{1*}\L.
\end{align*}
Notons que $\a,\b$ sont $U_I$-\'equivariants. En effet, on sait par d\'efinition que $\a$ est $U_I$-\'equivariant. Si $\L=\ZM/n$ pour un entier $n$ premier \`a $p,$ $\b$ est alors induit par le morphisme $\b(1):\L\to j_{1,I}^*R^1j_{1*}\L(1)$ qui envoie $1$ vers la classe du $\mu_n$-torseur des racines $n$-i\`emes de $t,$ o\`u $t$ est une \'equation locale de $\o{\pi}^{-1}_1(C_I).$ On peut se borner \`a un voisinage local strictement hens\'elien et $U_I$ stable. Un \'el\'ement $g\in U_I$ envoie ce torseur vers la classe du $\mu_n$-torseur des racines $n$-i\`emes de $t',$ o\`u $t'=tu$ avec $u$ une unit\'e. Dans ce voisinage l'\'equation $X^n-u=0$ admet une solution, donc le torseur donn\'e par $t'$ est isomorphe au torseur donn\'e par $t.$ C'est-\`a-dire $\b(1)$ est $U_I$-\'equivariant; il s'ensuit que $\b$ est $U_I$-\'equivariant. Donc, $e_1(j_{1,I}^*R^mj_{1*}\L)=\L$ pour $m=0,1.$ On a alors un triangle distingu\'e
$$
\L\To{}j_{1,I}^*Rj_{1*}\L\To{}\L(-1)[-1]\To{+1}
$$
dans $D^+(\o{\pi}^{-1}_1(C_I),\L-U_I),$ et on en d\'eduit que
$$
j_{1,I}^*Rj_{1*}\L=e_1(j_{1,I}^*Rj_{1*}\L)\      \ \text{et } \    \ e'_1(j_{1,I}^*Rj_{1*}\L)=0.
$$

Notons que le morphisme $\o{\rho}_{1,I}$ induit un foncteur $U_I$-\'equivariant:
$$\o{\rho}_{1,I,*}:D^+(\o{\pi}^{-1}_1(C_I),\L-U_I) \to D^+(Y^1_I\times\{0\},\L-U_I),$$ on a donc $e'_1(\o{\rho}_{1,I,*}j_{1,I}^*Rj_{1*}\L)=0.$

Par ailleurs,
\begin{align*}
\o{\rho}_{1,I,*}j_{1,I}^*Rj_{1*}\L&= i'^*_{1,I}Rj'_{1*}\rho_*\L\\
&= i'^*_{1,I}Rj'_{1*}(\L\oplus e'_1(\rho_*\L))\\
&=i'^*_{1,I}Rj'_{1*}\L\oplus i'^*_{1,I}Rj'_{1*}e'_1(\rho_*\L)\\
&=i'^*_{1,I}Rj'_{1*}\L\oplus e'_1(\o{\rho}_{1,I,*}j_{1,I}^*Rj_{1*}\L)\\
&=i'^*_{1,I}Rj'_{1*}\L,
\end{align*}
i.e. le morphisme (3) est un isomorphisme. D'o\`u l'\'enonc\'e du lemme.
\end{itemize}
\end{preuve}

\ali Notons que la composition $ \o{Y^1}\To{\o{\pi}^1} \o{Y^0}\into\DL^I$ induit un diagramme commutatif:
$$\xymatrix{
R\G(\DL,\L)^{U_I}\ar[r] \ar[d]_{(6)}^{\res.} &R\G(Y^0,\L)^{U_I}  \ar[r] \ar[d]_{(5)}^{\res.}  & R\G(Y^1,\L)^{U_I}\ar[d]^{\res.}_{\cong} \\
R\G(\o{\pi}^{-1}(C_I),j'^*_I Rj'_*\L)\ar[r] & R\G(\o{\pi}_0^{-1}(C_I),j^*_{0,I}Rj_{0 *}\L) \ar[r]  & R\G(\o{\pi}_1^{-1}(C_I),j^*_{1,I}Rj_{1 *}\L)
}$$

\begin{lem}
Les morphismes $(5)$ et $(6)$ sont des isomorphismes.
\end{lem}
\begin{preuve}
Consid\'erons tout d'abord le diagramme commutatif suivant:
$$
\xymatrix{
Y^1\ar@{^(->}[r]^{j_1}\ar[d]^{\pi^1} & \o{Y^1}\ar@{<-^)}[r]^{j_{1,I}}\ar[d]^{\o{\pi}^1} & \o{\pi}^{-1}_1(C_I)\ar[d]^{\o{\pi}^1_I}\\
Y^0 \ar@{^(->}[r]^{j_0} & \o{Y^0}\ar@{<-^)}[r]^{j_{0,I}} & \o{\pi}^{-1}_0(C_I)
}
$$
Par construction, $\o{\pi}^1:\o{Y^1}\to\o{Y^0}$ est un morphisme \'etale de groupe de Galois $N.$ Alors, pour tout faisceau \'etale $\FC$ sur $\o{Y^0},$ on a un isomorphisme canonique:
$$
R\G(\o{Y^0},\FC)\simto R_N(R\G(\o{Y^1},\o{\pi}^{1*}\FC)),
$$
o\`u $R_N$ est le foncteur d\'eriv\'e du foncteur des $N$-invariants.
En particulier, en tenant compte du fait que l'action de $N$ commute avec celle de $U_I,$ ceci induit un isomorphisme:
$$
R\G(Y^0,\L)^{U_I}\simto R_N(R\G(Y^1,\L)^{U_I}),
$$
 De m\^eme, $\o{\pi}^1_I:=\o{\pi}^1|_{\o{\pi}_1^{-1}(C_I)}$ induit un isomorphisme:
\begin{align*}
R\G(\o{\pi}_0^{-1}(C_I),j^*_{0,I}Rj_{0*}\L)&\simto R_N(R\G(\o{\pi}_1^{-1}(C_I),\o{\pi}^{1*}_I j^*_{0,I}Rj_{0*}\L))\\
&\simto R_N(R\G(\o{\pi}_1^{-1}(C_I), j^*_{1,I}Rj_{1*}\L)).
\end{align*}
Donc le morphisme (5) s'identifie \`a $R_N(R\G(Y^1,\L)^{U_I}\simto R\G(\o{\pi}_1^{-1}(C_I), j^*_{1,I}Rj_{1*}\L)).$
On en d\'eduit qu'il est un isomorphisme.

Notons que $\DL=Y^0\times_H T_d,$ et l'action de $U_I$ commute avec celle de $T_d.$ Alors
\begin{align*}
R\G(\DL,\L)^{U_I}
&=(R\G(Y^0,\L)\otimes^L_{\L[H]}\L T_d)^{U_I} \\
&=R\G(Y^0,\L)^{U_I}\otimes^L_{\L[H]}\L T_d
\end{align*}
D'autre part, on a $\L_{\DL}=\L_{Y^0}\otimes_{\L H}\L T_d.$ Alors,
$$
R\G(\o{\pi}^{-1}(C_I),j'^*_I Rj'_*\L)=R\G(\o{\pi}^{-1}_0(C_I),j^*_{0,I} Rj_{0,*}\L)\otimes^L_{\L[H]}\L T_d.
$$
On en d\'eduit que le morphisme (6) est un isomorphisme. Ceci termine la preuve du cas de codimension $1.$
\end{preuve}

\def\cprime{$'$} \def\cprime{$'$} \def\cprime{$'$} \def\cprime{$'$}
  \def\cprime{$'$}

\noindent
\textsc{Haoran Wang}\\
 Universit\'e Pierre et Marie Curie, Institut de Math\'ematiques de Jussieu

\noindent\texttt{haoran@math.jussieu.fr}

\medskip

\noindent{\it Adresse Pr\'esente:}
Max-Planck-Institut f\"ur Mathematik, Vivatsgasse 7, 53111 Bonn, Germany

\end{document}